\documentclass[12pt]{amsart}
\usepackage{amsfonts,amssymb}
\usepackage[all]{xypic}
\usepackage{amscd}

\newtheorem{theorem}{Theorem}[section]
\newtheorem{corollary}[theorem]{Corollary}
\newtheorem{Definition}[theorem]{Definition}
\newtheorem{lemma}[theorem]{Lemma}
\newtheorem{proposition}[theorem]{Proposition}
\newtheorem{Example}[theorem]{Example}
\newtheorem{Remark}[theorem]{Remark}

\newenvironment{remark}{\begin{Remark}\begin{em}}{\end{em}\end{Remark}}
\newenvironment{example}{\begin{Example}\begin{em}}{\end{em}\end{Example}}
\newenvironment{definition}{\begin{Definition}\begin{em}}{\end{em}\end{Definition}}

\numberwithin{equation}{section}

\newcommand{\R}{{\mathbb R}}

\newcommand{\N}{{\mathbb N}}
\newcommand{\ve}{\varepsilon}
\newcommand{\Bx}{{\mathbf x}}
\newcommand{\By}{{\mathbf y}}
\newcommand{\Bz}{{\mathbf z}}
\newcommand{\Bw}{{\mathbf w}}

\begin{document}

\title{A General Framework for Extending Means to Higher Orders}

\author{Jimmie Lawson and Yongdo Lim}

\address{Department of Mathematics, Louisiana State University,
Baton Rouge, LA70803, USA}

\email{lawson@math.lsu.edu}

\address{Department of Mathematics,
Kyungpook National University, Taegu 702-701, Korea}

\email{ylim@knu.ac.kr}

\date{\today}
\maketitle
\begin{abstract}In this paper we study the problem of extending means to
means of higher order.  We show how higher order means can be
inductively defined and established in general metric spaces, in
particular, in convex metric spaces.  As a particular application,
we consider the positive operators on a Hilbert space under
the Thompson metric and show that the operator logarithmic
mean admits  extensions of all higher orders, thus providing a
positive solution to a problem of Petz and Temesi \cite{PT}.
\end{abstract}

\section{Introduction}
A mean of order $n$ or $n$-mean for short on a set $X$ is a function $\mu:X^n\to X$
satisfying $\mu(x,x,\ldots, x)=x$
for all $x\in X$. Growing out of foundational papers on the subject such as that of Kubo and
Ando \cite{KA}, the theory of matrix and operator means of order $2$ has seen a
substantial development.  However,
no such developed theory has arisen for means of higher order.  Indeed there is no
 obvious definition for the extension of
a mean to the next higher order, although for various specific
means one can find a variety of candidates that have been put
forward.

Our purpose in this paper is to study a method of extending means to higher orders
that appears to offer a viable general approach.  For the definition and method of extending,
we favor  a general version of the recent approach of Horwitz \cite{Ho} for the case of means
 on the positive reals and of  Ando, Li, and Mathias \cite{ALM} for the case of the geometric mean on the
 positive (semi)definite Hermitian matrices.
  This approach has also been adopted and generalized beyond the geometric mean by Petz and
  Temesi in \cite{PT}, \cite{Pe}, although they only obtain existence for ordered tuples.
   In this paper we show that this approach can be generalized to means
on metric spaces and develop the theory of these extensions in this context.  The main theorems give
rather general conditions that guarantee that  extensions of all order exist.
Our main applications involve higher order means of positive operators on a Hilbert space.
  We are particularly interested in those cases in which one
starts with a mean of two variables and inductively extends it to all dimensions
greater than two.

In Section 2 we present our approach to extending means via limits of the ``barycentric operator.''
(This method is called ``symmetrization'' in \cite{PT}.)  Section 3 contains a major result of the paper:
means that are nonexpansive and coordinatewise contractive admit  extensions to all higher orders.
In Section 4 a special case of such means are considered, namely convex means for which the mean
is a metrically convex function assigning to any two points a metric midpoint. As explained in Section 5,
Hadamard spaces (for which the metric satisfies the semiparallelogram law) form an important class
of examples of metric spaces with associated convex mean.

In Section 6 we develop machinery for showing that certain types
of iterated means are nonexpansive and coordinatewise contractive.
Since many important means arise in this fashion, this is a useful
and important result.
Section 7 presents categorical aspects of means and their
extensions. In Section 8 a reverse construction is considered:
given a mean, when is it an extension of a lower order mean?
Connections between means, their extensions, and order are
developed in Section 9.  The paper closes in Sections 10 and 11
with the study of means on the space of positive operators on a
Hilbert  space. It is shown that certain iterated means are
nonexpansive and coordinatewise contractive, hence extend to
higher orders.  This is true, for example, of the logarithmic
mean, as we show in Section 11, and hence it has higher order
extensions, a conclusion that provides a positive solution of a
problem of Petz and Temesi in \cite{PT}.

\section{Mean Extensions}

An $n$-\emph{mean} on a set $X$ consists of an $n$-ary operation
(function) $\mu:X^n\to X$ that satisfies a generalized idempotency
law: $\mu(x,x,\ldots, x)=x$ for all $x\in X$.  A {\it mean} is
just an $n$-mean for some $n\geq 2$.  The mean is \emph{symmetric}
if it is invariant under permutations:
$$\mu(x_{\pi(1)},\ldots,x_{\pi(n)})=\mu(x_1,\ldots, x_n) \mbox{
for any permutation }\pi \mbox{ on } \{1,\ldots, n\}.$$ (Note that
for $2$-means this means that the binary operation given by the
mean is commutative.)  A \emph{topological $n$-mean} consists of a
Hausdorff topological space equipped with a continuous $n$-mean
$\mu$.  We also call the operation $\mu$ a topological $n$-mean.

A principal goal of our work is to ``extend''  a (symmetric)
$n$-mean on $X$ to a (symmetric) $n+1$-mean.  As mentioned in the
introduction we pattern our approach after \cite{Ho} and
\cite{ALM}, an approach that has some relation to the process of
compounding three means in three variables to obtain another such
mean (see \cite{BB}).

\begin{definition} \label{D:barycenter}
Given a set $X$ and a $k$-mean $\mu:X^k\to X$, the
\emph{barycentric operator} $\beta=\beta_\mu: X^{k+1}\to X^{k+1}$
is defined by
\begin{equation*} \beta(\Bx):=\bigl(\mu(\pi_{\ne 1}\Bx),\ldots,\mu(\pi_{\ne
k+1}\Bx)\bigr),
\end{equation*}
where $\Bx=(x_1,\ldots,x_{k+1})\in X^{k+1}$ and $\pi_{\ne j}
\Bx:=(x_1,\ldots,x_{j-1},x_{j+1},\ldots, x_{k+1})\in X^{k}$.

For a topological $k$-mean, we say that the barycentric map
$\beta$ is \emph{power convergent} if for each $\Bx\in X^{k+1}$,
we have $\lim_n \beta^n(\Bx)=(x^*,\ldots, x^*)$ for some $x^*\in
X$.
\end{definition}

As a motivating geometric example for the terminology
consider the $3$-mean in $\R^3$ that assigns to any three points
the centroid of the triangle for which they are the vertices, i.e., the point where the
three medians meet.  If we take now the four vertices of a $3$-simplex or tetrahedron
in $\R^3$,  the barycentric operator applied to the $4$-tuple consisting of the four vertices replaces
each vertex with the centroid (barycenter) of the face opposite it, the face with vertices the remaining
three vertices.  Thus one may envision the tetrahedron with vertices the four centroids of the four
faces as the result.  Repeating this process, one obtains a shrinking family of tetrahedra whose intersection
is the barycenter of the original tetrahedron, represented by the $4$-tuple with all entries equal to that point.

\begin{remark}\label{R:betastar}
There is an alternative way that the barycentric map of a $k$-mean may be defined, namely instead of defining the $i^{th}$-coordinate
of $\beta(\Bx)$ by deleting the $i^{th}$-coordinate of $\Bx\in X^{k+1}$ and applying $\mu$, we delete the coordinate
$i^*:=k+2-i$ and then apply $\mu$.    This means that we begin (from left to right) by first deleting coordinate $k+1$, $k$, down to $1$, instead
of beginning by deleting coordinate $1$ and continuing up to $k+1$.  We denote this alternative barycentric map by $\beta^*$.
One may define $\beta^*(\Bx)$ alternatively by reversing the $(k+1)$-tuple $\beta(\Bx)$.
Note, as long as the mean $\mu$ is symmetric, both methods
power  converge to the same limit, provided one of them power converges.
This equality of limit does not hold  in general for nonsymmetric means however.
The theories for $\beta$ and $\beta^*$ run parallel, so we restrict
our attention to $\beta$ with a few brief remarks concerning $\beta^*$.
\end{remark}

We define our first notion of an extension in terms of the
barycentric operator.

\begin{definition} \label{D:extend}
 A mean $\nu:X^{k+1}\to X$ is  a $\beta$-\emph{invariant extension of} $\mu: X^k\to X$ if
$\nu\circ\beta_\mu=\nu$, that is,
\begin{equation} \label{E:uniq}
\nu(\Bx)=\nu\bigl(\mu(\pi_{\ne 1}\Bx),\ldots,\mu(\pi_{\ne
{k+1}}\Bx)\bigr) \mbox{ for all } \Bx=(x_1,\ldots,x_{k+1})\in X^{k+1}.
\end{equation}

\end{definition}

The notion of a $\beta$-invariant extension was introduced by Horwitz \cite{Ho},
who called it type I invariance.

\begin{proposition}\label{P:unique}
Assume that $\mu:X^k\to X$ is a topological $k$-mean and that the
corresponding barycentric operator $\beta$ is power convergent.
Define $\tilde\mu:X^{k+1}\to X$ by $\tilde\mu(\Bx)=x^*$ where
$\lim_n\beta^n(\Bx)=(x^*,\ldots,x^*)$.
\begin{itemize}
\item[(i)]  $\tilde\mu:X^{k+1}\to X$ is a $(k+1)$-mean on $X$ that
is a $\beta$-invariant extension of $\mu$.
\item[(ii)] Any continuous mean on $X^{k+1}$ that is a $\beta$-invariant extension
of $\mu$ must equal $\tilde\mu$.
\item[(iii)] If $\mu$ is symmetric, so is $\tilde\mu$.
\end{itemize}
\end{proposition}

\begin{proof}  (i) For $x\in X$, $\Bx=(x,\ldots, x)\in X^{k+1}$,
we have $\beta(\Bx)=\Bx$ by the
idempotency of $\mu$.  Thus $\Bx=\lim_n \beta^n(\Bx)$ and hence
$\tilde\mu(\Bx)=x$, i.e., $\tilde\mu$ is a mean.  Further, we have
$$\tilde\mu(\beta(\Bx))=\pi_1\bigl(\lim_n\beta^n(\beta(\Bx))\bigr)
=\pi_1\bigl(\lim_n\beta^{n+1}(\Bx)\bigr)=
\tilde\mu(\Bx),$$ where $\pi_1$ is projection into the first
coordinate. Thus $\tilde\mu\circ\beta=\tilde\mu$, i.e.,
$\tilde\mu$ defines a $\beta$-invariant extension of $\mu$.

(ii) Suppose that $\nu$ is a continuous $(k+1)$-mean of $X$ that
is a $\beta$-invariant extension of $\mu$.  Since $\nu=\nu\circ
\beta=\nu\circ \beta^n$ (by repeated application of the first
equality),  for $\Bx\in X^{k+1},$
$$\nu(\Bx)=\nu(\beta(\Bx))=\nu(\beta^{n}(\Bx))=\nu(x^{*},\dots,
x^{*})=x^{*}=\tilde\mu(\Bx),$$ where
$(x^{*},\dots,x^{*})=\lim_{n\to\infty}\beta^{n}(\Bx).$

(iii) If $\mu$ is symmetric, then $\beta$
commutes with any permutation applied to the entries of $\Bx\in
X^{k+1}$, hence also $\beta^n$, and thus one obtains the same
limit with constant entry $\tilde\mu(x)$ in either case.  Hence
$\tilde\mu:X^{k+1}\to X$ is symmetric.
\end{proof}

\begin{remark}
If $\beta^*$ is power convergent, then a $\beta^*$-invariant extenstion is defined in the manner
of the previous proposition and the analogous proposition holds for $\beta^*$.
However, $\beta$-invariant extensions need not be $\beta^*$-invariant and
vice-versa.  But both notions collapse to the same one for symmetric means.
\end{remark}

We seek a notion of mean extension that allows one both to deduce
readily that a large number of properties transfer from a mean to
its extension and also is applicable to a wide variety of means.
The preceding proposition provides the ingredients for this
definition.

\begin{definition}\label{D:betaext}
A $(k+1)$-mean $\nu$ is a $\beta$-\emph{extension} of a topological
$k$-mean $\mu$ (or $\beta$-\emph{extends} $\mu$) if for each
$\Bx\in X^{k+1}$, $\lim_n \beta^n(\Bx)=(\nu(\Bx),\ldots,
\nu(\Bx))$.  In this case we say that $\beta$ \emph{power converges
to} $\nu$, written $\beta_\mu^n\to \nu$.
\end{definition}

We restate parts of Proposition \ref{P:unique} in terms of this
definition.

\begin{corollary}\label{C:betaext}
If $\beta_\mu$ power converges, where $\mu$ is a topological mean,
then it converges to a  $(k+1)$-mean $\tilde\mu$,
which $($by definition$)$ is a $\beta$-extension of $\mu$.
Furthermore, $\tilde\mu$, if continuous, is the unique $\beta$-invariant extension of
$\mu$.
\end{corollary}

\begin{remark} \label{R:horwitz}
A.\ Horwitz \cite{Ho} and later D.\ Petz and R.\ Temesi \cite{PT} consider means on the positive reals
and show that any continuous symmetric $2$-mean that is strict
($\min(a,b)<\mu(a,b)<\max(a,b)$ for $a\ne b$) and order-preserving
in each variable has a power convergent barycentric map, and hence
has a unique $\beta$-extension to a $3$-mean. Petz and Temesi point out
that the argument for power convergence extends to higher order variables,
and thus one can inductively define $\beta$-extensions for all $n>2$
\cite[Section 5]{PT}.
For the arithmetic,
geometric, and harmonic means the extensions yield the usual
corresponding means of $n$-variables.  To check this,
one has only to note that they are continuous and are
$\beta$-invariant extensions, then apply the previous corollary.
\end{remark}

\section{Power Convergence}

In this section we consider properties preserved by
$\beta$-extensions and develop sufficient conditions for a
topological mean to (recursively) admit a $\beta$-extension.

Given $X$ equipped with a $k$-mean $\mu$, a subset $C$ is
\emph{convex} if $\mu(x_1,\ldots,x_k)\in C$ whenever
$x_1,\ldots,x_k\in C$.

\begin{lemma}\label{L:convex}
If a topological mean admits a $\beta$-extension, then any closed
set that is convex with respect to the mean is convex with respect
to the extension.
\end{lemma}

\begin{proof} Let $\mu:X^k\to X$ be the given mean, and let
$x_1,\ldots,x_{k+1}\in A$, a closed $\mu$-convex set.  Set
$\Bx:=(x_1,\ldots,x_{k+1})$. Then by convexity each coordinate of
$\beta(\Bx)$ is in $A$ and by induction each coordinate of
$\beta^n(\Bx)$ is in $A$.  Since $A$ is closed, it follows that
the coordinate limits, which are all
$\tilde\mu(x_1,\ldots,x_{k+1})$, belong to $A$, where $\tilde\mu$
is the $\beta$-extension.
\end{proof}

Recall that the \emph{convex hull} of a set $A$ is the smallest
convex subset containing $A$, and can be obtained by intersecting
all convex sets containing $A$. In a similar fashion in the case
of a topological mean $\mu$ the \emph{closed convex hull} can be
obtained by intersecting all closed convex sets containing $A$ or,
as follows from the continuity of $\mu$, by closing up the convex
hull.

\begin{definition}\label{D:loccon}
A topological mean is \emph{locally convex} if there exists at
each point a basis of (not necessarily open) neighborhoods that
are convex.  A metric topological mean is \emph{uniformly locally
convex} if for each $\ve>0$, there exists $\delta>0$ such that the
diameter of the convex hull of $A$ is less than $\ve$ whenever the
diameter of $A$ is less than $\delta$.  A metric topological mean
is \emph{closed ball convex} if all closed balls
${\overline{B}}_\ve(x):=\{y\in X:d(x,y)\leq \ve\}$ are convex for
all $x\in X$.
\end{definition}

\begin{remark}
Since in a topological mean the closure of a set $A$ is convex
whenever $A$ is and has the same diameter, the convex
sets in the definition of locally convex and uniformly locally
convex may be taken to be closed if the space is regular, which is
the case in the metric setting.
\end{remark}

\begin{lemma}\label{L:hier}
Given a topological mean on a metric space, closed ball convexity
implies uniform local convexity, which in turn implies local
convexity.
\end{lemma}

\begin{proof} Let $\ve>0$.  Choose $\delta=\ve/4$.  Then any set
$A$ of diameter less than $\delta$ is contained in a closed ball
of radius less than $2\delta$ around any point of $A$, which in
turn has diameter less than $4\delta=\ve$. If $X$ is closed ball
convex, then this closed ball is convex and hence contains the
closed convex hull of $A$.    The proof that uniformly locally
convex implies locally convex is straightforward.
\end{proof}

 The next lemma is an immediate consequence of Lemma
\ref{L:convex}.

\begin{lemma}\label{L:loccon}
Let $X$ be a topological mean with a $\beta$-extension. If $X$ is
locally convex resp.\ metric and uniformly locally convex resp.\
metric and closed ball convex with respect to the given mean, then
it is with respect to the extension.
\end{lemma}

\begin{lemma}\label{L:context}
If $\nu$ is a $\beta$-extension of the topological mean $\mu$ and
if $X$ is locally convex and regular, then $\nu$ is continuous.
\end{lemma}

\begin{proof}
Let $\mu:X^k\to X$ and let $\nu:X^{k+1}\to X$ be the
$\beta$-extension.  Let $\Bx=(x_1,\ldots,x_{k+1})\in X^{k+1}$, let
$x^*=\nu(\Bx)$, and let $U$ be an open set containing $x^*$.  Pick
a closed convex neighborhood $V$ of $x^*$ such that $V\subseteq
U$. Since by hypothesis the sequence $\beta^n(\Bx)$ converges to
the $k+1$-string with entries $x^*$, we have $\beta^n(\Bx)\in
V^{k+1}$ for some $n$ large enough.  By continuity of $\mu$ and
hence of $\beta^n$, there exists $W$ open in $X^{k+1}$ containing
$\Bx$ such that $\beta^n(W)\subseteq V^{k+1}$.  For any $\By\in
W$, we have $\beta^n(\By)\in V^{k+1}$, and hence $\beta^m(\By)\in
V^{k+1}$ for all $m>n$ since $V$ is convex.  Since $V$ is closed
it follows that $\nu(\By)\in V$.  Thus $\nu$ is continuous.
\end{proof}

\begin{definition} Let $\mu:X^k\to X$ be $k$-mean on a metric space $X$.
For $\Bx=(x_1,\ldots,x_{k+1})\in X^{k+1}$, set
$\vert\Bx\vert=\{x_1,\ldots,x_{k+1}\}$, the underlying set of the
$(k+1)$-tuple, and define the \emph{diameter} $\Delta(\Bx)$ of
$\Bx$ by
$$\Delta(\Bx)=\mbox{diam}\vert \Bx\vert=\sup\{d(x_i,x_j):1\leq
i,j\leq k+1\}.$$ The mean $\mu$ is \emph{weakly
$\beta$-contractive} if for each $\Bx\in X^{k+1}$, we have $\lim_n
\Delta(\beta^n(\Bx))=0$.  For $0<\rho<1$, we say that $\mu$ is
\emph{coordinatewise} $\rho$-\emph{contractive} if  for any
$\Bx,\By\in X^k$ that differ only in one coordinate, say $x_j\ne
y_j$,
$$d(\mu(\Bx),\mu(\By))\leq\rho\, d(x_j,y_j)$$
\end{definition}

\begin{lemma}\label{L:rctowc}
If $\mu:X^k\to X$ is a coordinatewise $\rho$-contractive mean for
$0<\rho<1$, then it is weakly $\beta$-contractive.
\end{lemma}

\begin{proof} Assume that $\mu:X^k\to X$ is a coordinatewise
$\rho$-contractive $k$-mean. We equip $X^{k+1}$ with the sup
metric
$$d\big((x_1,\ldots,x_{k+1}),(y_1,\ldots,y_{k+1})\bigr):=\max\{d(x_j,y_j):
1\leq j\leq k+1\}.$$  We show by induction on $n$ that for any
$\Bx\in X^{k+1}$
 and any two adjacent coordinates $(\beta^n(\Bx))_i$ and
  $(\beta^n(\Bx))_{i+1}$, $$d( (\beta^n(\Bx))_i,(\beta^n(\Bx))_{i+1})
  \leq \rho^n d(x_i,x_{i+1}).$$

For $n=1$, we note by the coordinatewise $\rho$-contractive
property that
\[
d( (\beta(\Bx))_i,(\beta(\Bx))_{i+1})= d\bigl(\mu(\pi_{\ne
i}(\Bx)),\mu(\pi_{\ne i+1} (\Bx)\bigr) \leq \rho d(x_i,x_{i+1}),\]
since $\pi_{\ne i}(\Bx)$ and $\pi_{\ne i+1}(\Bx)$ differ only in
the $i^{th}$-coordinate, where they have the entries $x_{i+1}$ and
$x_i$ resp.  Assume the validity of the inductive hypothesis for
$n$. Since $\beta^{n+1}(\Bx)=\beta(\beta^n(\Bx))$, we have from
the case $n=1$ that
$$d\bigl( (\beta^{n+1}(\Bx))_i,(\beta^{n+1}(\Bx))_{i+1}\bigr)
\leq \rho\cdot d\bigl((\beta^n(\Bx))_i,(\beta^n(\Bx)
)_{i+1}\bigr).$$ By the inductive hypothesis, the latter is less
than or equal $\rho\cdot\rho^n
d(x_i,x_{i+1})=\rho^{n+1}d(x_i,x_{i+1})$. This completes the
induction.  We then conclude from the triangle inequality that
since any two entries of $\beta^n(\Bx)$ are at most $k$ steps
apart, we have $d(\beta^n(\Bx)_i,\beta^n(\Bx)_j)\leq k\rho^n
\Delta(\Bx)$, and thus $\Delta(\beta^n(\Bx))\leq
k\rho^n\Delta(\Bx)$. Therefore $\lim_n \Delta(\beta^n(\Bx))=0$.
\end{proof}

Note that if $\beta_\mu$ is power convergent, then $\mu$ must be
weakly $\beta$-contractive. The next proposition provides a
converse.

\begin{proposition}\label{P:contract}
Let $X$ be a complete metric space endowed with a weakly
$\beta$-contractive $k$-mean $\mu$.  If $X$ is uniformly locally
convex, then $\beta$ is power convergent, so that a
$\beta$-extension exists.
\end{proposition}

\begin{proof}  For $\Bx\in X$, set $C_n(\Bx)$ equal to the
closed convex hull of $\vert\beta^n(\Bx)\vert$.  By hypothesis
$\Delta(\beta^n(\Bx))=\mbox{diam}\vert \beta^n(\Bx)\vert\to 0$ and
then by uniform local convexity\ \, diam $C_n(\Bx)\to 0$. Note
that since $C_n(\Bx)$ is convex, it contains $\vert
\beta^m(\Bx)\vert$ for all $m>n$, and hence contains $C_m(\Bx)$.
Thus the collection $\{ C_n(\Bx)\}$ is a decreasing sequence of
closed convex sets whose diameters converge to $0$.  Since $X$ is
a complete metric space the intersection consists of a single
point $\{x^*\}$, and it is now easy to show that $\beta^n(\Bx)$
converges to the $(k+1)$-tuple with all entries $x^*$.
\end{proof}

We next single out another property that will be useful for inductively
$\beta$-extending a mean and give some useful equivalences.

\begin{lemma}\label{L:lip1}
Let $X$ be a metric space endowed with a $k$-mean $\mu$. We endow
$X^k$ and $X^{k+1}$ with the sup metric that takes the supremum of
the distances between each of the corresponding coordinates.  Then
the following three conditions are equivalent:
\begin{itemize}
\item[(1)] for all $\Bx=(x_1,\ldots,x_k),\By=(y_1,\ldots,y_k)\in X^k$,
\begin{equation*}
d(\mu(\Bx),\mu(\By))\leq \max\{d(x_j,y_j):1\leq j\leq k\};
\end{equation*}
\item[(2)] The mean $\mu:X^k\to X$ is Lipschitz with Lipschitz
constant $1\ ($hence, in particular, is continuous).
\item[(3)] The map $\beta: X^{k+1}\to X^{k+1}$ is Lipschitz with
Lipschitz constant $1$.
\end{itemize}
These conditions imply
\begin{itemize}
\item[(4)] $X$ is closed ball convex.
\end{itemize}
\end{lemma}

\begin{proof} $(1)\Leftrightarrow(2)$: the right-hand side of
(1) is the definition of the sup metric, so the two statements are
equivalent.

$(2)\Rightarrow (3)$: In each coordinate the map $\beta$ is a
projection followed by $\mu$, a composition of maps with Lipschitz
constant $1$, and thus has Lipschitz constant $1$.  Since this
holds in each coordinate, it holds in the sup metric.

$(3)\Rightarrow (2)$:
$(3)\Rightarrow (2)$:  Fixing some $z\in X$, we have for
$\Bx\in X^k$ that
$\mu(\Bx)=\pi_1(\beta(z,\Bx))$, and the right-hand side is a
composition of maps of Lipschitz constant $1$.

$(2)\Rightarrow (4)$: For $\ve>0$ and $x\in X$, $y_1,\ldots,y_k\in
X^k$, we have for $\By:=(y_1,\ldots,y_k)$
\[d(x,\mu(\By))\leq d\bigl( \mu(x,\ldots,x),\mu(\By)\bigr)\leq
d\bigl( (x,\ldots,x),\By)=\max_i d(x,y_i)\leq \ve\] provided
$d(x,y_i)\leq \ve$ for all $i$.  Thus $X$ is closed ball convex.
\end{proof}

\begin{definition}\label{D:nonexp}
A $k$-mean $\mu$ on a metric space $X$ is called
\emph{nonexpansive} if it satisfies for all
$\Bx=(x_1,\ldots,x_k),\By=(y_1,\ldots,y_k)\in X^k$,
\begin{equation}
d(\mu(\Bx),\mu(\By))\leq \max\{d(x_j,y_j):1\leq j\leq k\},
\end{equation}
or equivalently condition (2) or (3) of the preceding lemma.
\end{definition}

\begin{lemma}\label{L:nonexp}
If $\mu$ is a nonexpansive $k$-mean on a metric space $X$ and if
$\mu$ has a $\beta$-extension $\tilde\mu$, then $\tilde\mu$ is
nonexpansive.
\end{lemma}

\begin{proof} Let $\pi_1:X^{k+1}\to X$ denote projection into the
first coordinate.  For $\Bx\in X^{k+1}$,
$$\tilde\mu(\Bx)=\pi_1(\lim_n \beta^n(\Bx))=\lim_n
(\pi_1\circ\beta^n)(\Bx).$$ Since $\tilde\mu$ is the pointwise
limit of Lipschitz maps $\pi_1\circ\beta^n$ of Lipschitz constant
one (by Lemma \ref{L:lip1}), it is also.
\end{proof}

The next proposition is the principal tool that allows us to
extend means inductively to higher order.

\begin{proposition}\label{P:extend}
Let $X$ be a complete metric space equipped with a nonexpansive,
coordinatewise $\rho$-contractive $(0<\rho<1)$\ $k$-mean
$\mu:X^k\to X$, $k\geq 2$.  Then the barycentric operator $\beta$
is power convergent, and hence there exists a (unique) continuous
$(k+1)$-mean $\tilde\mu:X^{k+1}\to X$ that $\beta$-extends $\mu$.
Furthermore, $\tilde\mu:X^{k+1}\to X$ is nonexpansive and
coordinatewise $\rho$-contractive.
\end{proposition}

\begin{proof}  By Lemma \ref{L:rctowc} $\mu$ is weakly contractive.
Since $\mu$ is nonexpansive, by Lemma \ref{L:lip1} $X$ is closed
ball convex, hence uniformly locally convex (Lemma \ref{L:hier}),
and thus $\beta$ is power convergent and has a $\beta$-extension
to a $(k+1)$-mean $\tilde\mu$ by Proposition \ref{P:contract}.  By
Lemma \ref{L:nonexp} the mean $\tilde\mu$ is nonexpansive and
hence continuous (Lemma \ref{L:lip1}(2)).

To finish we show that $\tilde\mu$ is coordinatewise
$\rho$-contractive.  Let $\Bx,\By\in X^{k+1}$ differ only in the
$j^{th}$-coordinate, $x_j\ne y_j$. Then by coordinatewise
$\rho$-contractivity
$$d\bigl((\beta(\Bx))_i,(\beta(\By))_i\bigr)=d\bigl(\mu(\pi_{\ne
i}(\Bx)),\mu(\pi_{\ne i}(\By))\bigr)\leq \rho d(x_j,y_j),$$ since
$\pi_{\ne i}(\Bx)$ and $\pi_{\ne i}(\By)$ differ in at most one
coordinate, and are then $x_j$ and $y_j$ in that coordinate. Since
the inequality holds for each $i$, we have
$d(\beta(\Bx),\beta(\By))\leq \rho d(x_j,y_j)$. Since $\beta$ is
nonexpansive by Lemma \ref{L:lip1}, we conclude that
$d(\beta^n(\Bx),\beta^n(\By))\leq \rho d(x_j,y_j)$ for all $n$.
Taking the limit as $n\to \infty$, we obtain $d(\tilde\mu(\Bx),
\tilde\mu(\By))\leq \rho d(x_j,y_j)$.
\end{proof}

The next theorem is the culmination of this section.  It follows
from a straightforward induction using the preceding proposition.

\begin{theorem}\label{T:big}
Let $X$ be a complete metric space equipped with a nonexpansive,
coordinatewise $\rho$-contractive $(0<\rho<1)$\ $k$-mean
$\mu:X^k\to X$, $k\geq 2$.  Then there exists uniquely a family of
continuous means $\mu_n:X^n\to X$, one for every $n>k$, such that
each is a $\beta$-extension of the previous one.  Furthermore,
each $\mu_n$ is nonexpansive and coordinatewise
$\rho$-contractive.
\end{theorem}

\begin{example}
Consider on $\R$ the mean $\mu(x,y)=sx+(1-s)y$, where $0<s<1$.  Set $\rho=\mbox{max}\{s,1-s\}$.
Then it is an elementary calculation to verify that $\mu$ is coordinatewise $\rho$-contractive  and
nonexpansive.  Hence $\mu$ inductively $\beta$-extends to an $n$-mean for all $n>2$.  For example,
if $m(x,y)=(2/3)x+(1/3)y$, then one varifies that $m(x,y,z)=(2/5)x+(7/20)y+(1/4)z$
is a $\beta$-invariant extension of $m$, and hence must be its three-variable $\beta$-extension.
If one uses the alternative barycentric operator $\beta^*$, then one obtains the extension
$m_3(x,y,z)=(4/7)x+(2/7)y+(1/7)z$.
\end{example}



\section{Convex means}

In general a metric space may have none, one, or many midpoints
between two given points in the space.  (Recall that $m$ is a
midpoint of $a$ and $b$ if $d(m,a)=d(m,b)=(1/2)d(a,b)$.) We wish
to consider the setting where possibly many midpoints may exist,
but there is a distinguished midpoint, and these distinguished
midpoints appear in a ``convex'' manner.

\begin{definition}
A symmetric mean $\mu:X\times X\to X$, written $\mu(x,y)=x\#y$, on
a complete metric space $X$ is called a \emph{convex mean} if it
satisfies the \emph{basic convexity condition}
\begin{eqnarray}\label{bcc}
d(x\#z,y\#z)\leq \frac{1}{2}d(x,y) \mbox{ for all }x,y,z\in X.
\end{eqnarray}
\end{definition}

\begin{lemma}\label{L:midpt}
For a convex mean, $x\#y$ is a metric midpoint for all $x,y$.
\end{lemma}

\begin{proof}
By the basic convexity condition $d(x\#y,y=y\#y)\leq (1/2)d(x,y)$
and similarly $d(x\#y,x)\leq (1/2)d(x,y)$. Thus
$$d(x,y)\leq d(x,x\#y)+d(x\#y,y)\leq
\frac{1}{2}d(x,y)+\frac{1}{2}d(x,y)=d(x,y).$$ It follows that each
inequality is an equality, so $d(x,x\#y)+d(x\#y,y)=d(x,y)$. Hence
adding together the two inequalities in the first line of the
proof gives an equality, so each inequality is an equality.
\end{proof}

Note that in the case there is only one metric midpoint between
two points $x,y$ it must be that case that $x\#y$ is that
midpoint.

The next proposition gives a useful equivalence for convexity.
\begin{proposition}\label{P:equiv1}
Let $(X,d)$ be a complete metric space equipped with a mean $\mu$.
Then $\mu$ is a convex mean if and only if
$$d(x\#y,u\#v)\leq \frac{1}{2}d(x,u)+\frac{1}{2}d(y,v) \mbox{ for
all } x,y,u,v\in X.$$
\end{proposition}

\begin{proof}
For a convex mean $d(x\#y,u\#y)\leq (1/2)d(x,u)$ and
$d(x\#y,u\#v)\leq (1/2)d(y,v)$.  The condition of the theorem now
follows by adding the equations and an application of the triangle
inequality.

The reverse implication follows by choosing, $x=x$, $z=y=v$, and
$y=u$.
\end{proof}

\begin{proposition}\label{P:conv3}
A convex mean inductively $\beta$-extends to a symmetric,
nonexpansive, coordinatewise $(1/2)$-contractive $n$-mean for
every $n>2$.
\end{proposition}

\begin{proof}
Note that the definition of a convex mean is that of a symmetric
coordinatewise $(1/2)$-contractive $2$-mean.  Proposition
\ref{P:equiv1} further yields that it is nonexpansive, since
$$d(u\#y,u\#v)\leq \frac{1}{2}d(x,u)+\frac{1}{2}d(y,v)\leq
\max\{d(x,u),d(y,v)\}.$$ Thus by Theorem \ref{T:big} we obtain
inductively a $\beta$-extension for every $n$ that is nonexpansive
and coordinatewise $(1/2)$-contractive. By Proposition
\ref{P:unique} each extension is symmetric.
\end{proof}

\begin{example}\label{E:banach}
Let $X$ be a Banach space (or a closed convex subset thereof) and define the symmetric
 $2$-mean $\mu(x,y)=(1/2)(x+y)$.
This is the midpoint with respect to the norm metric, and is
easily seen to be a convex mean.  Setting
$\mu_k(x_1,\ldots,x_k)=(1/k)\sum_{i=1}^k x_i$,  one verifies
directly that $\mu_{k+1}$ is the $\beta$-extension of $\mu_k$, so
$\mu$ inductively $\beta$-extends to the standard arithmetic mean
$\mu_n$ for all $n$.
\end{example}

\section{Hadamard spaces}

A metric space $X$ is said to satisfy the \emph{semiparallelogram law} if for any two points
$x_1,x_2,\in X$, there exists $z\in X$ that satisfies for all $x\in X$:
$$d(x_1,x_2)^2+4d(x,z)^2\leq 2d(x,x_1)^2+2d(x,x_2)^2.$$
It follows readily that $z$ is the unique midpoint between $x_1$
and $x_2$.  A \emph{Hadamard space} (occasionally called a
\emph{Bruhat-Tits space}) is a complete metric space that
satisfies the semiparallelogram  law.

Using a metric notion for an upper bound of curvature (geodesic
triangles in the metric space satisfy certain inequalities when
compared with test triangles), one calls a metric space a
CAT($\kappa$)-\emph{space} if it is a geodesic space (each pair of
points can be connected by a metric geodesic) satisfying the
curvature bound condition for the real number $\kappa$  (see
\cite[Chapter i]{Ball} or \cite[Chapter II.1]{BrH}). The
CAT($0$)-spaces are the \emph{non-positively curved} spaces.  A
metric space has an alternative characterization as a Hadamard
space: it is a simply connected, complete,  geodesic
CAT($0$)-space (see \cite[Proposition 5.1, Chapter 1]{Ball} or
\cite[Exercise 1.9, Chapter II.1]{BrH}).

What is important for our current purposes is the following remark:

\begin{remark}\label{R:hadamard}
Let $X$ be a Hadamard metric space and define a $2$-mean by $\mu(x,y)$ is the unique
midpoint between $x$ and $y$.  This defines a convex mean in the sense of the preceding
section, see \cite[Proposition 5.4, Chapter I]{Ball} or \cite[Proposition 2.2, Chapter II.2]{BrH}.  Hence
by the preceding section this mean may be $\beta$-extended to an $n$-mean for every $n>2$.
\end{remark}

A wide variety of Hadamard spaces and constructions for new Hadamard spaces from old appear
in \cite{Ball} and \cite{BrH}.   Some examples include Hadamard manifolds (simply connected complete
Riemannian manifolds with nonpositive sectional curvature), particularly simply connected symmetric
spaces of noncompact type, finite-dimensional hyperbolic geometries over the
reals, complexes, and quaternions, symmetric cones, Tits buildings, and various examples
obtained by coning and gluing.

Of particular interest to us is the example of the manifold of
positive definite matrices endowed with the usual Riemannian
metric called the trace metric.  This metric yields a Hadamard
manifold and the midpoint mean operation in this case is precisely
the geometric mean of the two positive definite matrices; see
\cite{LL} and the references there. Using the fact that the length
metric satisfies the semiparallelogram law, hence is a convex
metric with the midpoint operation being a convex mean, we obtain
the following alternative derivation of the principal result of
\cite{ALM}:

\begin{corollary}\label{C:trace}
Let $X$ denote the set of positive definite real or complex
matrices equipped with the Riemannian trace metric.  Then the
midpoint operation for the corresponding length metric, which is
precisely the geometric mean, defines a convex $2$-mean, which
$($by Proposition \ref{P:conv3}$)$ $\beta$-extends to an $n$-mean
for each $n>2$.
\end{corollary}

\section{Iterated means}

A standard construction technique for means is iteration, the
arithmetic-geometric mean being the best known example.  In this
section we develop machinery for showing that certain iterated
means are coordinatewise $\rho$-contractive and nonexpansive,
hence admit $\beta$-extensions of all orders. We apply this
machinery to operator means in a later section.

\begin{definition}\label{D:iteratedmean}
Let $\lambda,\nu$ be $2$-means on a complete metric space $X$.
Starting with $\lambda_1=\lambda$ and $\nu_1=\nu$, we define
inductively the sequences of means $\{\lambda_n\}$ and $\{\nu_n\}$
by
\begin{eqnarray*} \lambda_{n+1}(x,y)&=&\lambda(\lambda_n(x,y),\nu_n(x,y)),
\quad \nu_{n+1}(x,y)=\nu(\lambda_n(x,y),\nu_n(x,y))\\
(resp.\ \nu_{n+1}(x,y) &=& \nu(\lambda_n(x,y),\nu_n(x,y)),\quad
\lambda_{n+1}(x,y)=\lambda(\lambda_n(x,y),\nu_{n+1}(x,y)).
\end{eqnarray*}  If there exists a $2$-mean $\mu$ such that for all $x,y\in
X$, $\lim_n \lambda_n(x,y)=\mu(x,y)=\lim_n \nu_n(x,y)$, then $\mu$
is called the \emph{iterated composition} (resp.\ \emph{skewed
iterated composition}) of $\lambda$ and $\nu$ and denoted
$\mu=\lambda*\nu$ (resp.\ $\mu=\lambda*_s\nu$).
\end{definition}

We begin with a useful lemma that ensures convergence.
\begin{lemma}\label{P:convergence}
Let $\{x_n\},\{y_n\}$ be sequences in a complete metric space $X$
satisfying one of the two following conditions:
\begin{itemize}
\item[(i)] for each $k\geq 1$,  $x_{k+1}$ is a midpoint of $x_k$ and $y_k$ and $d(x_{k+1},y_{k+1})\leq d(x_{k+1},y_k)$;  or
\item[(ii)] for each $k\geq 1$,  $x_{k+1}$ is a midpoint of $x_k$ and $y_{k+1}$
and $d(x_k,y_{k+1})\leq d(x_k,y_k)$.
\end{itemize}
Then both sequences are Cauchy and converge to the same point.
\end{lemma}

\begin{proof}
Assume (i).  For any $n\geq 1$, we have by hypothesis
$$ d(x_{n+1},y_{n+1})\leq d(x_{n+1},y_n)=\frac{1}{2}d(x_n,y_n),$$
where the last equality follows from that the fact that $x_{n+1}$
is a midpoint of $x_n,y_n$.  Similarly
$$d(x_{n+1},x_n)=\frac{1}{2}d(x_n,y_n).$$
It follows by induction resp.\ induction and the triangle
inequality that
$$d(x_{n},y_n)\leq \frac{1}{2^{n-1}}d(x_1,y_1)\mbox{ resp. }
d(x_{n+k},x_n)\leq \biggl(\sum_{i=0}^{k-1}\frac{1}{2^{n+i}} \biggr)
d(x_1,y_1)<\frac{1}{2^{n-1}}d(x_1,y_1).$$ Thus the sequence
$\{x_n\}$ is Cauchy, and hence converges, and the sequence $\{y_n\}$
must also approach the same limit.

Part (ii) follows by applying part (i) to the sequences $\{x_n\}$
and $\{z_n\}$, where $z_n=y_{n+1}$.

\end{proof}

\begin{proposition}\label{P:iterconv}
Let $\lambda$ be a convex mean and $\nu$ be a nonexpansive mean on
a complete metric space $X$.  Then the iterated composition $\mu=
\lambda*\nu$ resp.\ the skewed iterated composition
$\mu=\lambda*_s \nu$ exists and is nonexpansive.
\end{proposition}

\begin{proof}
For $x,y\in X$, we set $x_n=\lambda_n(x,y)$ and $y_n=\nu_n(x,y)$
(see Definition \ref{D:iteratedmean}). Then
$x_{n+1}=\lambda(x_n,y_n)$ and $y_{n+1}=\nu(x_n,y_n)$.  We observe
that
$$d(x_{n+1},y_{n+1})=d(\nu(x_{n+1},x_{n+1}),\nu(x_{n},y_{n}))
\leq \max\{d(x_{n+1},x_n),d(x_{n+1},y_n )\},$$ where the last
inequality follows from the fact that $\nu$ is nonexpansive. Since
$\lambda$ is a convex mean $x_{n+1}=\lambda(x_n,y_n)$ is a
midpoint for $x_n$ and $y_n$, hence
$d(x_{n+1},x_n)=d(x_{n+1},y_n)$, and thus $d(x_{n+1},y_{n+1})\leq
d(x_{n+1},y_n)$, i.e., condition (i) of Lemma \ref{P:convergence}
is satisfied.  It thus follows that $\lim_n x_n=\lim_n y_n$
exists, and we define this limit to be $\mu(x,y)$.  If $x=y$, then
it is immediate that $x=x_n=y_n$ for all $n$, so $\mu$ is a mean.
Thus the iterated composition $\mu=\lambda*\nu$ exists.

For the case of the skewed iterated mean, we set
$x_1=\lambda(x,y)$, $y_1=\nu(x,y)$ and $$y_{n+1}=\nu(x_n,y_n),
\quad x_{n+1}=\lambda(x_n,y_{n+1}).$$ Then $$d(x_k,y_{k+1})=
d(\nu(x_k,x_k),\nu(x_k,y_k))\leq \max\{d(x_k,x_k),d(x_k,y_k)\}
=d(x_k,y_k),$$ where the inequality follows from the nonexpansive
property.  Thus Lemma \ref{P:convergence}(ii) is satisfied. That
the skewed iterated composition $\mu=\lambda*_s \nu$ exists now
follows as in the preceding paragraph.

It follows from Proposition \ref{P:equiv1} that $\lambda$ is
Lipschitz with Lipschitz constant $1$ and the same holds for $\nu$
since it is non-expansive.  In both cases the higher numbered
means $\lambda_n$ and $\nu_n$ are built up from these by products
and compositions, so are also $1$-Lipschitz (recall that product
metrics are always the sup metric). Since $\mu$ is the pointwise
limit of the sequence $\{\lambda_n\}$ (and $\{\nu_n\}$), it is
also $1$-Lipschitz, i.e., nonexpansive.
\end{proof}

\begin{proposition}\label{P:approximate}
Suppose in a complete metric space $X$ that $\lambda$ is a convex
mean and $\nu$ is coordinatewise $\rho'$-contractive, $0<
\rho'<1,$ and nonexpansive.  Then the iterated composition
$\lambda*\nu$ resp.\ the skewed iterated composition
$\lambda*_s\nu$ exists and is coordinatewise $\rho$-contractive,
$\rho=\max\{1/2,\rho'\}$, and nonexpansive, hence $\beta$-extends
to all orders greater than two.
\end{proposition}

\begin{proof}
By Proposition \ref{P:iterconv} the iterated composition and
skewed iterated composition both exist and are nonexpansive.

We establish that $\mu=\lambda*\nu$ is coordinatewise
$\rho$-contractive. Let $a,b,c\in X$. To calculate $\mu(a,b)$ and
$\mu(a,c)$ we define inductively
\begin{eqnarray*} b_1^-&=&\lambda(a,b),\
b_1^+=\nu(a,b),\ b_{k+1}^-=\lambda(b_k^-,b_k^+),\ b_{k+1}^+=\nu(b_k^-,b_k^+)\\
c_1^-&=&\lambda(a,c),\ c_1^+=\nu(a,c),\
c_{k+1}^-=\lambda(c_k^-,c_k^+),\ c_{k+1}^+=\nu(c_k^-,c_k^+).
\end{eqnarray*}
Note that $b_k^-=\lambda_k(a,b)$, $b_k^+=\nu_k(a,b)$,
$c_k^-=\lambda_k(a,c)$, $c_k^+=\nu_k(a,c)$. We have
$d(b_1^-,c_1^-)=d(\lambda(a,b),\lambda(a,c))\leq (1/2)d(b,c)$ by
convexity of $\lambda$ and similarly $d(b_1^+,c_1^+)\leq \rho
d(b,c)$ by coordinatewise $\rho$-contractivity of $\nu$.

We claim that by induction $$d(b_n^-,c_n^-)\leq
 \rho d(b,c)
\mbox{ and } d(b_n^+,c_n^+)\leq \rho d(b,c).$$  By the preceding
paragraph it holds for $n=1$.  Assume that it is true for $n=k$.
Then
\begin{eqnarray*} d(b_{k+1}^-,c_{k+1}^-) &=& d(\lambda(b_k^-,b_k^+),\lambda(c_k^-,
c_k^+))\\
&\leq &d(\lambda(b_k^-,b_k^+),\lambda(c_k^-,
b_k^+))+d(\lambda(c_k^-,b_k^+),\lambda(c_k^-,c_k^+))\\
 &\leq & \frac{1}{2}d(b_k^-,c_k^-)+\frac{1}{2}d(b_k^+,c_k^+)\\
&\leq & \frac{1}{2}(\rho d(b,c)+\rho d(b,c))=\rho d(b,c).
\end{eqnarray*}
 Using the
nonexpansivity of $\nu_{k+1}$, we obtain
 \begin{eqnarray*} d(b_{k+1}^+,c_{k+1}^+)&=& d(\nu(b_k^-,b_k^+),
\nu(c_k^-,c_k^+))\\
& \leq & \mbox{max}\{ d(b_k^-,c_k^-),d(b_k^+,c_k^+)\}\\
& \leq & \mbox{max}\{ \rho d(b,c), \rho d(b,c)\}=\rho d(b,c).
\end{eqnarray*}
This completes the induction.  Note that in the alternative
notation we have shown that $d(\lambda_n(a,b),\lambda_n(a,c))\leq
\rho d(b,c)$ and $d(\nu_n(a,b),\nu_n(a,c))\leq \rho d(b,c)$ for
all $n\in \N$.

Since $\mu=\lambda*\nu$,  $\lim_n b_n^-=\lim_n
\lambda_n(a,b)=\mu(a,b)$, $\lim_n c_n^-=\lim_n
\lambda_n(a,c)=\mu(a,c)$. By continuity of $d(\cdot,\cdot)$ and
the preceding paragraph, it follows that $d(\mu(a,b),\mu(a,c)))
\leq \rho d(b,c)$.

The proof that the skewed iterated composition $\mu=\lambda*_s\nu$
is coordinatewise $\rho$-contractive is similar, but contains a
twist or two. To calculate $\mu(a,b)$ and $\mu(a,c)$ for $a,b,c\in
X$, we define inductively
\begin{eqnarray*}
b_1^+&=&\nu(a,b),\ b_1^-=\lambda(a,b),\
b_{k+1}^+=\nu(b_k^-,b_k^+),
\ b_{k+1}^-=\lambda(b_k^-,b_{k+1}^+)\\
c_1^+&=&\nu(a,c),\ c_1^-=\lambda(a,c),\
c_{k+1}^+=\nu(c_k^-,c_k^+),\ c_{k+1}^-=\lambda(c_k^-,c_{k+1}^+).
\end{eqnarray*}
We have $d(b_1^-,c_1^-)=d(\lambda(a,b),\lambda(a,c))\leq
(1/2)d(b,c)$ by convexity of $\lambda$ and similarly
$d(b_1^+,c_1^+)\leq \rho d(b,c)$ by coordinatewise
$\rho$-contractivity of $\nu$.

We claim that by induction $$d(b_n^-,c_n^-)\leq \rho d(b,c) \mbox{
and } d(b_n^+,c_n^+)\leq \rho d(b,c).$$  By the preceding paragraph
it holds for $n=1$.  Assume that it is true for $n=k$. Using the
nonexpansivity of $\nu$, we obtain
 \begin{eqnarray*} d(b_{k+1}^+,c_{k+1}^+)&=& d(\nu(b_k^-,b_k^+),\nu(c_k^-,c_k^+))\\
& \leq & \mbox{max}\{ d(b_k^-,c_k^-),d(b_k^+,c_k^+)\}\\
& \leq & \mbox{max}\{ \rho d(b,c), \rho d(b,c)\}=\rho d(b,c).
\end{eqnarray*}
It then follows that
\begin{eqnarray*} d(b_{k+1}^-,c_{k+1}^-) &=& d(\lambda(b_{k}^-,b_{k+1}^+),\lambda(c_{k}^-,
c_{k+1}^+))\\
&\leq &d(\lambda(b_k^-,b_{k+1}^+),\lambda(c_k^-,
b_{k+1}^+))+d(\lambda(c_k^-,b_{k+1}^+),\lambda(c_k^-,c_{k+1}^+))\\
 &\leq & \frac{1}{2}d(b_k^-,c_k^-)+\frac{1}{2}d(b_{k+1}^+,c_{k+1}^+)\\
&\leq & \frac{1}{2}(\rho d(b,c)+\rho d(b,c))=\rho d(b,c).
\end{eqnarray*}
This completes the induction.

By hypothesis $\lim_n b_n^-=\lim_n \lambda_n(a,b)=\mu(a,b)$,
$\lim_n c_n^- =\lim_n \lambda_n(a,c) =\mu(a,c)$.  By continuity of
$d(\cdot,\cdot)$ and the preceding paragraph, it follows that
$d(\mu(a,b),\mu(a,c))) \leq \rho d(b,c)$.

The last assertion of the the proposition now follows from Theorem
\ref{T:big}.

\end{proof}

\section{Categorical constructions}

In this section we consider the behavior of mean extensions with
respect to standard constructions such as continuous images,
products, and subspaces.

\begin{definition}A function $g:(X,\mu)\to (Y,\nu),$ where $\mu, \nu$ are
$k$-means on $X$ and $Y$ respectively, is called a \emph{$k$-mean
homomorphism} or \emph{homomorphism} for short,  if $g\circ
\mu=\nu\circ g_{k}$, that is, the following diagram commutes
\[
\begin{CD}
X  @>g>> Y \\
@A{\mu}AA  @AA{\nu}A\\
X^{k}  @>{g_{k}}>>  Y^{k} \\
\end{CD}
\]
where $g_{k}:X^{k}\to Y^{k},\
g_{k}(x_{1},\dots,x_{k}):=(g(x_{1}),\dots,g(x_{k})).$
\end{definition}

\begin{proposition}\label{P:homo}
Let $(X,\mu)$ and $(Y,\nu)$ be topological $k$-means,  and let $g:X\to Y$ be
a continuous $k$-mean homomorphism.

\begin{itemize}
\item[(i)] If each of $\mu$ and $\nu$ $\beta$-extend to
$(k+1)$-means $\tilde \mu$ and\ $\tilde\nu$ resp.,   then $g:X\to Y$ is a
$(k+1)$-mean homomorphism.
\item[(ii)] If $g$ is surjective and $\mu$ $\beta$-extends to a $(k+1)$-mean $\tilde\mu$, then
$\nu$ $\beta$-extends to a $(k+1)$-mean $\tilde\nu$, and $g$ is then a $(k+1)$-mean
homomorphism.
\end{itemize}
\end{proposition}

\begin{proof} (i) It follows directly from the fact that $g$ is $k$-mean homomorphism that
$\beta_Yg_{k+1}=g_{k+1}\beta_X:X^{k+1}\to Y^{k+1}$ (indeed commutation of $g_k$ with $\beta$
is an equivalence).  By induction $$\beta_Y^ng_{k+1}=g_{k+1}\beta_X^n:X^{k+1}\to Y^{k+1}$$  for
all $n>0$.  Taking the limit of both sides as $n\to\infty$ and projecting into the first coordinate yields
(i).

(ii) For $\By\in Y^{k+1}$, there exists $\Bx\in X^{k+1}$ such that $g_{k+1}(\Bx)=\By$.
Again we have $$\beta_Y^n(\By)=\beta_Y^ng_{k+1}(\Bx)=g_{k+1}\beta_X^n(\Bx)$$ for $n>0$.  By hypothesis
the right-hand side converges to a diagonal element with entries $g(\tilde\mu(\Bx))$ as $n\to\infty$, so that the left-hand side
also converges to a diagonal element.  Thus $\nu$ $\beta$-extends to $\tilde\nu$.  The last assertion follows from (i).
\end{proof}

\begin{example}\label{E:harmonic} Let
$$\mu_f(x,y)=f^{-1}\bigg(\frac{f(x)+f(y)}{2}\bigg)$$
be a \emph{quasi-arithmetic mean} defined on $\R^+$ by a continuous strictly monotonic function $f$ \cite{PT}.
By the preceding proposition applied to $g=f^{-1}$ and Example \ref{E:banach}
$$\mu_k(x_1,\ldots,x_k)=f^{-1}\big(\frac{1}{k}\sum_{i=1}^k f(x_i)\big).$$
Note that the arithmetic, geometric, and harmonic means belong to
this class by taking the identity map, the logarithmic map, and
the inversion map (on the positive reals) respectively.  More
generally one can take on the positive reals the generalized or
power mean $m(x,y)=\big((x^\alpha+y^\alpha)2\big)^{1/\alpha}$ for
$\alpha\ne 0$ with $f(x)=x^\alpha$.

An analogous construction and characterization of the higher order means remains valid for the power means
on the space of positive definite matrices.  Note that the case $\alpha=1$ gives the arithmetic mean and the case
$\alpha=-1$ gives the harmonic mean.
\end{example}

\begin{definition}
Let $\mu, \nu$ be  $k$-means on $X$ and $Y,$ respectively. Define
$\mu\times \nu:(X\times Y)^{k}\to X\times Y$  by
$$(\mu\times\nu)(x_{1},y_{1},x_{2},y_{2},\dots,
x_{k},y_{k})=(\mu(x_{1},x_{2},\dots,x_{k}),\nu(y_{1},y_{2},\dots,y_{k})).
$$ Then $\mu\times \nu$ is a $k$-mean on $X\times Y,$ called the
\emph{product mean} of $\mu$ and $\nu.$ Indeed, if
$\Bz=(x,y,x,y,\dots,x,y)\in (X\times Y)^{k}$ for fixed $(x,y)\in
X\times Y,$ then $(\mu\times
\nu)(\Bz)=(\mu(x,x,\dots,x),\nu(y,y,\dots,y))=(x,y).$
\end{definition}

\begin{theorem}\label{T:product}
Let $(X,d_{1})$ and $(Y,d_{2})$ be complete metric spaces equipped
with nonexpansive, coordinatewise $\rho$-contractive resp.\
$\rho'$-contractive $(0<\rho, \rho'<1)$\ $k$-means $\mu:X^k\to X$
resp.\  $\nu:Y^{k}\to Y$. Then the mean $\mu\times \nu$
is a nonexpansive, coordinatewise $\max\{\rho,
\rho'\}$-contractive $k$-mean on $X\times Y.$ Furthermore, for
$n\geq k,$  its $n$-mean $\beta$-extension $(\mu\times \nu)_{n}$ coincides
with the product mean $\mu_{n}\times \nu_{n}$ of the individual $\beta$- extension
$n$-means:
$$(\mu\times \nu)_{n}=\mu_{n}\times\nu_{n}.$$
\end{theorem}

\begin{proof}
For $\Bz=(x_{1},y_{1},x_{2},y_{2},\dots, x_{k},y_{k})\in (X\times
Y)^{k},$ we denote $\Bz_{x}=(x_{1},x_{2},\dots, x_{k})\in X^{k}, \
\ \Bz_{y}=(y_{1},y_{2},\dots, y_{k})\in Y^{k}.$ Then $(\mu\times
\nu)(\Bz)=(\mu(\Bz_{x}), \nu(\Bz_{y})).$

 Let $\Bz=(x_{11},y_{11},x_{12},
y_{12},\dots, x_{1k}, y_{1k}), \Bw=(x_{21},y_{21},x_{22},
y_{22},\dots, x_{2k}, y_{2k})\in (X\times Y)^{k}.$ Then by
nonexpansive property of $\mu$ and $\nu,$
\begin{eqnarray*}
d\Big((\mu\times \nu)(\Bz),(\mu\times
\nu)(\Bw)\Big)&=&d\Big((\mu(\Bz_{x}), \nu(\Bz_{y})),(\mu(\Bw_{x}),
\nu(\Bw_{y}))\Big)\\
&=&\max\{d_{1}(\mu(\Bz_{x}),\mu(\Bw_{x})),
d_{2}(\nu(\Bz_{y}),\nu(\Bw_{y}))\}\\
&\leq&\max\{\max\{d_{1}(x_{1j},x_{2j})\},
\max\{d_{2}(y_{1j},y_{2j})\}:1\leq j\leq k\}\\
&=& \max\{d_{1}(x_{1j},x_{2j}),
d_{2}(y_{1j},y_{2j}):1\leq j\leq k\}\\
&=&\max \{d((x_{1j},y_{1j}), (x_{2j},y_{2j})):1\leq j\leq k\},
\end{eqnarray*} which implies that $\mu\times \nu$ is a non-expansive
$k$-mean on $X\times Y$ equipped with the sup-metric.

If $\Bz$ and $\Bw$  differ only in  one coordinate of  $(X\times
Y)^{k},$ say $(x_{1j},y_{1j})\neq (x_{2j},y_{2j})$, but
$(x_{1i},y_{1i})= (x_{2i},y_{2i})$, $1\leq i\neq j\leq k,$ then the
inequality in preceding argument turns into
\begin{eqnarray*}
d\Big((\mu\times \nu)(\Bz),(\mu\times
\nu)(\Bw)\Big)&=&d\Big((\mu(\Bz_{x}), \nu(\Bz_{y})),(\mu(\Bw_{x}),
\nu(\Bw_{y}))\Big)\\
&=&\max\{d_{1}(\mu(\Bz_{x}),\mu(\Bw_{x})),
d_{2}(\nu(\Bz_{y}),\nu(\Bw_{y}))\}\\
&\leq&\max\{\rho d_{1}(x_{1j},x_{2j}),
\rho'd_{2}(y_{1j},y_{2j})\}\\
&\leq&\max\{\rho,\rho'\}\max\{ d_{1}(x_{1j},x_{2j}),
d_{2}(y_{1j},y_{2j})\}\\
&\leq&\max\{\rho,\rho'\}d\Big((x_{1j},y_{1j}),(x_{2j},y_{2j})\Big).
\end{eqnarray*}
Therefore $\mu\times \nu$ is a coordinatewise $\max\{\rho,\rho'\}$-
contractive $k$-mean on $X\times Y.$

Next, we will prove $(\mu\times \nu)_{n}(\Bz)=(\mu_{n}(\Bz_{x}),
\nu(\Bz_{y})),\Bz\in (X\times Y)^{n}, n\geq k$ by induction. The
case $n=k$ follows  by the definition of $\mu\times \nu.$ Suppose
that the assertion  holds true for $n-1.$ Let $\gamma:(X\times
Y)^{n}\to X\times Y$ be defined by $\gamma(\Bz)=(\mu_{n}(\Bz_{x}),
\nu_{n}(\Bz_{y})).$ Then $\gamma$ is continuous and hence it
suffices to show $\gamma\circ \beta_{n}=\gamma$ by the uniqueness
of mean extension (Proposition \ref{P:unique}) where $\beta_{n}$
is the barycentric operator on $(X\times Y)^{n}$ obtained from the
$(n-1)$-mean $(\mu\times \nu)_{n-1}.$

For $\Bz=(z_{1},z_{2},\dots,z_{n})\in (X\times Y)^{n},
z_{i}=(x_{i},y_{i})\in X\times Y,$ we have from induction
$$(\mu\times \nu)_{n-1}(\pi_{\neq i}\Bz)=
(\mu_{n-1}(\pi_{\neq i}\Bz_{x}),\nu_{n-1}(\pi_{\neq i}\Bz_{y})), \
1\leq i\leq n,$$ and then from $\beta_{n}(\Bz)=\Big((\mu\times
\nu)_{n-1}(\pi_{\neq 1}\Bz),(\mu\times \nu)_{n-1}(\pi_{\neq
2}\Bz), \dots, (\mu\times \nu)_{n-1}(\pi_{\neq n}\Bz)\Big),$
\begin{eqnarray*}
\beta_{n}(\Bz)_{x}&=&\Big(\mu_{n-1}(\pi_{\neq 1}\Bz_{x}),
\mu_{n-1}(\pi_{\neq
2}\Bz_{x}),\dots, \mu_{n-1}(\pi_{\neq n}\Bz_{x})\Big)=\beta_{\mu}(\Bz_{x})\\
\beta_{n}(\Bz)_{y}&=&\Big(\nu_{n-1}(\pi_{\neq 1}\Bz_{y}),
\nu_{n-1}(\pi_{\neq
2}\Bz_{y}),\dots, \nu_{n-1}(\pi_{\neq n}\Bz_{y})\Big)=\beta_{\nu}(\Bz_{y})\\
\end{eqnarray*}
where $\beta_{\mu}, \beta_{\nu}$ denote the barycentric operators
of the $n$-means $\mu_{n}, \nu_{n}$; $\mu_{n}\circ
\beta_{\mu}=\mu_{n}$ and $ \nu_{n}\circ \beta_{\nu}=\nu_{n}.$ Now,
\begin{eqnarray*}
 (\gamma\circ
 \beta_{n})(\Bz)&=&\gamma(\beta_{n}(\Bz))=\Big(\mu_{n}(\beta_{n}(\Bz)_{x}),
 \nu_{n}(\beta_{n}(\Bz)_{y})\Big)\\
 &=&\Big(\mu_{n}(\beta_{\mu}(\Bz_{x})),
 \nu_{n}(\beta_{\nu}(\Bz_{y}))\Big)\\
 &=&(\mu_{n}(\Bz_{x}), \nu_{n}(\Bz_{y}))=\gamma(\Bz)
 \end{eqnarray*}
 which completes the claim.
 \end{proof}

 \begin{remark} The product mean satisfies the associative law
 $$\mu\times (\nu\times \omega)=(\mu\times \nu)\times \omega$$ for
 any  $k$-means $\mu, \nu,$ and $\omega$ on
 $X,Y,Z.$ If these are nonexpansive and coordinatewise contractive,
 then its $n$-mean extension satisfies
$$(\mu\times\nu\times \omega)_{n}=\mu_{n}\times \nu_{n}\times \omega_{n}.$$
\end{remark}

The next result of the section is quite straightforward and hence the proof is omitted.

\begin{proposition}\label{P:subspace}
If $Z$ is a nonempty closed $k$-submean of a topological $k$-mean
$(X,\mu)$ that $\beta$-extends, then $Z$ also $\beta$-extends and
$\tilde\mu\vert Z^{k+1}= \widetilde{\mu\vert_Z^k}$.
\end{proposition}

\begin{corollary}\label{C:local}
Suppose that a $k$-mean $(X,\mu)$ restricts to a nonexpansive, coordinatewise
$\rho_n$-contractive mean on $A_n$, where $0<\rho_n<1$ for each $n$ and
$A_n$ is an increasing sequence of closed convex sets with $X=\bigcup_n A_n$.
Then $\mu$ inductively $\beta$-extends to $j$-mean for each $j>k$ in such a way
that the restriction to $A_n$ is the appropriate $\beta$-extension of the restriction
of $\mu$ to $A_n$.
\end{corollary}

\begin{proof}  Given any $j$-tuple in $X^{j}$ for $j>k$, there exists some $A_m$ that
contains all entries of the tuple.  Applying Theorem \ref{T:big} to the restriction of $\mu$ to $A_m$,
we conclude that $\mu\vert_{A_m}$ $\beta$-extends for every index greater than $k$.  Thus in
particular the appropriate extension exists to evaluate the given $j$-tuple.  It is clear
that if a larger $A_{m+n}$ is chosen, one obtains the same calculaton.  Thus the $\beta$-extension
is independent of the containing $A_m$, and hence we obtain a $\beta$-extension of $\mu$ on all of $X$.
\end{proof}

\section{Stability and reductions}

\begin{definition} A $k$-mean $\mu$ on a set $X$ is called $\beta$-\emph{stable}
if the graph of $\mu$ is invariant under $\beta$, that is,
$$\forall \Bx\in X^{k+1},\pi_{k+1}(\Bx)=\mu(\pi_{\ne  k+1}(\Bx))\Rightarrow
\pi_{k+1}(\beta(\Bx))=\mu(\pi_{\ne k+1}(\beta(\Bx))).$$
A $(k+1)$-mean $\nu$ on $X$ is called a \emph{stable extension} of $\mu$ if
$$\nu(x_1,\ldots, x_k,\mu(x_1,\ldots,x_k))=\mu(x_1,\ldots,x_k)\mbox{ for all
}x_1,\ldots,x_k
\in X.$$ Conversely $\mu$ is called a \emph{stable reduction} of $\nu$.
\end{definition}

In \cite{Ho} Horwitz says that $\nu$ is type 2 invariant with respect to $\mu$ if
$\nu$ is a stable extension of $\mu$.  He only considers the case of a $2$-mean $\mu$
and a $3$-mean $\nu$.

\begin{proposition}\label{P:stable}
If a topological $k$-mean $\mu$ is $\beta$-stable and admits a $\beta$-extension
$\tilde\mu$, then $\tilde\mu$ is a stable extension of $\mu$.
\end{proposition}

\begin{proof}
Note that  $\beta$ does not change the last coordinate of any
$k+1$-tuple $\Bx=(x_1,\ldots,x_k,\mu(x_1,\ldots,x_k))$ in the
graph of $\mu$. Hence if $\mu$ is $\beta$-stable, it follows that
$\beta^n(\Bx)$ has the same last coordinate for all $n$. Since
$\beta$ power converges to $\tilde\mu$, it follows that the
diagonal limiting value has entries the last coordinate of $\Bx$,
namely $\mu(x_1,\ldots,x_k)$.  Thus
$\tilde\mu(\Bx)=\mu(x_1,\ldots,x_k)$, so $\tilde\mu$ is a stable
extension of $\mu$.
\end{proof}

\begin{proposition}\label{P:reduction}
If $X$ is a complete metric space equipped with a
coordinatewise $\rho$-contractive $($symmetric$)$ $(k+1)$-mean
$\nu$, $k\geq 2$ then $\nu$ admits exactly one stable reduction
$($which is also symmetric$)$.
\end{proposition}

\begin{proof}
Suppose that  $X$ is a complete metric space equipped with a
coordinatewise $\rho$-contractive $(k+1)$-mean $\nu$. The map
$g:X\to X$ defined by $g(x)=\nu(x_1,\ldots,x_{k},x)$ is by
hypothesis $\rho$-contractive, and hence has a unique fixed point.
Define $\mu(x_1,\ldots,x_k)$ to be this fixed point. It follows
immediately that $\mu$ is a stable reduction of $\nu$, and
uniqueness of the fixed point guarantees that it is the unique
reduction.  Note that the fact that $\nu$ is a mean and the
definition of $\mu$ imply that $\mu$ is also a mean. The symmetry
of $\mu$ follows directly from that of $\nu$.
\end{proof}

The following corollary is an immediate consequence of
Propositions \ref{P:extend}, \ref{P:stable}, and
\ref{P:reduction}.

\begin{corollary}\label{C:reduction}
Let $X$ be a complete metric space equipped with a nonexpansive, coordinatewise $\rho$-contractive
$k$-mean $\mu$.  If $\mu$ is $\beta$-stable, then the unique $\beta$-extension $\tilde\mu$ is
a stable extension of $\mu$, and $\mu$ is the unique stable reduction of $\tilde\mu$.
\end{corollary}

 \begin{definition}
A $2$-mean $\mu$ on a set $X$ satisfies the \emph{limited medial
property} if $\mu(a,b)=\mu(x,y)=:m$ implies that
$\mu(\mu(a,x),\mu(b,y))=m.$ For $m\in X,$
$$X_{m}:=\{(x,y)\in X^{2}: \mu(x,y)=m\}.$$
\end{definition}

\begin{remark}\label{R:limmed}
(1) If a $2$-mean $\mu(x,y)=x\#y$ satisfies the limited medial
property, then it is $\beta$-stable since
$$\beta(x,y,x\#y)=(y\#(x\#y),x\#(x\#y),x\#y).$$  The latter is in
the graph of $\mu$, since by limited mediality $x\#y=(x\#y)\#(x\#
y)$ implies $(x\#(x\#y))\#(y\#(x\# y))=x\# y$.

\par\noindent
(2) It was shown in \cite{ALM} that the matrix geometric mean for
positive definite matrices satisfies the limited medial property.
This was extended to very general notions of geometric mean in
\cite{LL5}, in particular for the geometric mean of positive
operators on a Hilbert space or, more generally, for the positive
elements of a $C^*$-algebra.  Hence by part (1) and the earlier
results each $k$-extension of the geometric mean yields both the
higher order ones by stable extension and the lower order ones by
stable reduction.
\end{remark}

\begin{lemma} The mean $\mu$  satisfies the limited medial property
if and only if $(\mu\times \mu)(X_{m}\times X_{m})\subset X_{m}$ for
each $m\in X.$
\end{lemma}

\begin{proof}  Let $(a,b),(x,y)\in X_{m},$ that is $\mu(a,b)=\mu(x,y)=m.$
Then $(\mu(a,x),\mu(b,y))\in X_{m}$ (limited medial property) if and
only if $\mu(\mu(a,x),\mu(b,y))=(\mu\times \mu)((a,b),(x,y))=m.$
\end{proof}

\begin{proposition}Let $X$ be complete metric space equipped with
nonexpansive, coordinatewise $\rho$-contractive $2$-mean
$\mu(x,y)=x\#y$ satisfying the limited medial property.  If
$x_{i}\#y_{i}=m$ for all $1\leq i\leq n,$ then
$\mu_{n}(\Bx)\#\mu_{n}(\By)=m,$ where
$\Bx=(x_{1},\dots,x_{n}),\By=(y_{1},\dots,y_{n})\in X^{n}.$
\end{proposition}

\begin{proof} By previous lemma $(\mu\times \mu)(X_{m}\times X_{m})\subset X_{m},$ and
$X_{m}$ is a non-empty closed subset of $X^{2}.$ Let
$\omega:=(\mu\times\mu)\Big|_{X_{m}^{2}}.$ By Theorem
\ref{T:product}(3), the $n$-mean extension $\omega_{n}:X_{m}^{n}\to
X_{m}$ of $\omega$ is given by
$$\omega_{n}=(\mu_{n}\times\mu_{n})\Big|_{X_{m}^{n}}.$$
Suppose that $x_{i}\#y_{i}=m, i=1,2,\dots,n.$ Then
$\Bz:=(x_{1},y_{1},x_{2},y_{2},\dots,x_{n},y_{n})\in X_{m}^{n}$ and
so $\omega_{n}(\Bz)=(\mu_{n}(\Bx),\mu_{n}(\By))\in X_{m}$ which
implies that $\mu_{n}(\Bx)\#\mu_{n}(\By)=m.$
\end{proof}

\begin{remark}Under the assumption of the preceding proposition, if
$a\#b=x\#y=m$ then $X_{m}\subseteq \{(z,w)\in X\times X:
\mu_{3}(a,b,z)\#\mu_{3}(x,y,w)=m\}.$
\end{remark}

\section{Ordered convex metric spaces}
Throughout this section we assume that $X$ is a complete metric
space equipped with a non-expansive, coordinatewise
$\rho$-contractive $k$-mean $\mu:X^{k}\to X, k\geq 2,$ and $\mu_{n}$
($n>k$) denotes the non-expansive, coordinatewise $\rho$-contractive
 $n$-mean obtained inductively.  We further assume that $X$ is equipped a closed partial order $\leq.$
We let $\leq_{k}$ be the product order on $X^{k}$ defined by
$$(x_{1}, \dots, x_{k})\leq_{k} (y_{1}, \dots, y_{k})\ \mathrm{if\ and\ only\
if\ }\  x_{i}\leq y_{i}, 1\leq i\leq k.$$  Recall  for a map $g:X\to Y,$
we let $g_{k}:X^{k}\to Y^{k}, (x_{1},\dots, x_{k})\mapsto
(g(x_{1}),\dots, g(x_{k})).$

\begin{definition} A $k$-mean $\nu$ on $X$ is  said to be \emph{monotone}
for the partial order order
$\leq$ if $\nu(\Bx)\leq \nu(\By)$ for any $\Bx, \By\in X^{k }$ with
$\Bx\leq \By.$
\end{definition}

\begin{theorem}\label{T:mono} (1)
If a nonexpansive, coordinatewise $\rho$-contractive $k$-mean $\mu$
 is  monotone for the closed partial order $\leq,$
then $\mu_{n}$ is also monotone for any  $n\geq k$:
$$\Bx\leq_{n} \By\Longrightarrow \mu_{n}(\Bx)\leq \mu_{n}(\By).$$
(2)  Let $(X,\delta)$ be another complete metric space equipped
with a nonexpansive, coordinatewise $\rho$-contractive $k$-mean
$\nu:X^{k}\to X.$ Suppose that  $\leq$ is closed in the product
topology $(X,d)\times (X,\delta)$  and   either $\mu$ or $\nu$ is
monotone with respect to  $\leq.$ Then $\mu\leq \nu$ implies $
\mu_{n}\leq \nu_{n}$ for any $n\geq k.$

\end{theorem}

\begin{proof}
(1) We induct on $m$ beginning at $k$.  Suppose that $\mu_{m}$  is monotone. Let $\Bx,
\By\in X^{m+1}$ with $\Bx\leq_{m+1}\By.$ Then $\pi_{\neq
i}\Bx\leq_{m}\pi_{\neq i}\By$ for all $1\leq i\leq m+1.$ By the
induction hypothesis, $\mu_{m}(\pi_{\neq i}\Bx)\leq\mu_{m}(\pi_{\neq
i}\By)$  for each $i,$ and hence
$\beta_{m+1}(\Bx)\leq_{m+1}\beta_{m+1}(\By).$  By repeated application, we have
$\beta_{m+1}^{n}(\Bx)\leq_{m+1}\beta_{m+1}^{n}(\By),\
n=1,2,\dots.$ By the closeness of the order, $\lim_{n\to
\infty}\beta_{n+1}^{n}(\Bx)\leq_{n+1}\lim_{n\to\infty}\beta_{n+1}^{n}(\By).$
In particular, $\mu_{n+1}(\Bx)\leq \mu_{n+1}(\By).$

(2)  Assume that $\mu_{m}(\Bx)\leq \nu_{m}(\Bx)$ for all $\Bx\in
X^{m}$ for some  $m\geq k.$ Let $\Bx\in X^{m+1}.$ Then $\mu_{m}(\pi_{\neq
j}\Bx)\leq \nu_{m}(\pi_{\neq j}\Bx)$ for $1\leq j\leq m+1$ and hence
\begin{eqnarray}\label{E:1}  \beta(\Bx)\leq_{m+1}
\alpha(\Bx) \end{eqnarray}where $\beta$ and $\alpha$ are the
barycentric operators with respect to the  means $\mu$ and $\nu,$
respectively. We will show by induction that
$\beta^{n}(\Bx)\leq_{m+1} \alpha^{n}(\Bx)$ for all $n.$  Suppose
that it holds true for $n.$ Then $\pi_{\neq j}\beta^{n}(\Bx)\leq
\pi_{\neq j}\alpha^{n}(\Bx)$ for all $1\leq j\leq m+1.$
If $\nu$ is monotone, $\nu_{m}$ is monotone by (1) and
\begin{eqnarray}\label{E:2}
\nu_{m}\bigl(\pi_{\neq j}\beta^{n}(\Bx)\bigr)\leq
\nu_{m}\bigl(\pi_{\neq j}\alpha^{n}(\Bx)\bigr), \ \ 1\leq j\leq m+1.
\end{eqnarray}
Therefore,
\begin{eqnarray*}
\beta^{n+1}(\Bx)&=& \beta(\beta^{n}(\Bx))\stackrel{\eqref{E:1}}\leq_{m+1} \alpha(\beta^{n}(\Bx))\\
&=&\bigl(\nu_{m}(\pi_{\neq 1}\beta^{n}(\Bx)), \dots,
\nu_{m}(\pi_{\neq
m+1}\beta^{n}(\Bx))\bigr)\\
&\stackrel{\eqref{E:2}}\leq_{m+1}&\bigl(\nu_{m}(\pi_{\neq
1}\alpha^{n}(\Bx)), \dots, \nu_{m}(\pi_{\neq
k+1}\alpha^{n}(\Bx))\bigr)\\
&=&\alpha(\alpha^{n}(\Bx))=\alpha^{n+1}(\Bx).
\end{eqnarray*}

In the case  that $\mu$ is monotone, (\ref{E:2})
changes to
\begin{eqnarray}\label{E:3}
\mu_{m}\bigl(\pi_{\neq j}\beta^{n}(\Bx)\bigr)\leq
\mu_{m}\bigl(\pi_{\neq j}\alpha^{n}(\Bx)\bigr), \ \ 1\leq j\leq k+1,
\end{eqnarray} and $\beta^{n+1}(\Bx)\leq_{m+1}
\alpha^{n+1}(\Bx)$ from
\begin{eqnarray*}
\beta^{n+1}(\Bx)&=& \beta(\beta^{n}(\Bx))\\
&=&\bigl(\mu_{m}(\pi_{\neq 1}\beta^{n}(\Bx)), \dots,
\mu_{m}(\pi_{\neq m+1}\beta^{n}(\Bx))\bigr)\\
&\stackrel{\eqref{E:3}}\leq_{m+1}&\bigl(\mu_{m}(\pi_{\neq
1}\alpha^{n}(\Bx)), \dots,
\mu_{m}(\pi_{\neq m+1}\alpha^{n}(\Bx))\bigr)\\
&\stackrel{{\mathrm{induction}}}\leq_{m+1}& \bigl(\nu_{m}(\pi_{\neq
1}\alpha^{n}(\Bx)), \dots, \nu_{m}(\pi_{\neq
m+1}\alpha^{n}(\Bx))\bigr)\\
&=&\alpha(\alpha^{n}(\Bx))=\alpha^{n+1}(\Bx).
\end{eqnarray*}

Since the order is closed in the product topology, the
 inequality holds for limits:
$$\lim_{n\to\infty}\beta^{n}(\Bx)\leq_{m+1}
 \lim_{n\to \infty}\alpha^{n}(\Bx).$$ In particular,
 $\mu_{m+1}(\Bx)\leq \nu_{m+1}(\Bx).$
\end{proof}

\begin{theorem}\label{T:ineq} Let $X$ and $Y$ be complete metric spaces equipped
with  nonexpansive, coordinatewise $\rho$ (resp.\
$\rho'$)-contractive $k$-means $\mu$ and $\nu,$ respectively. Let
$\leq$ be a closed partial order on $Y.$ Suppose that the mean $\nu$
is monotone for the partial order $\leq$ on $Y.$ If $g:X\to Y$ is
continuous and satisfies $g\circ \mu\leq \nu\circ g_{k}$ $($resp.
$g\circ \mu\geq \nu\circ g_{k}),$ then for any $n\geq k,$ $$g\circ
\mu_{n}\leq \nu_{n}\circ g_{n} \ \ \ \bigl({\mathrm{resp.}}\ \
g\circ\mu_{n}\geq \nu_{n}\circ g_{n}\bigr).$$
\end{theorem}

\begin{proof}Suppose that $g\circ \mu\leq \nu\circ g_{k}$; the other
case is similar. Suppose that
$g\circ \mu_{m}\leq \nu_{m}\circ g_{m}$ holds true for some
$m\geq k.$ Let $\alpha$ be the barycentric operator with respect to
the mean $\nu$ on $Y.$ Then $g\bigl(\mu_{m}(\pi_{\neq
j}\Bx)\bigr)\leq \nu_{m}g_{m}(\pi_{\neq j}\Bx)$ for $\Bx\in X^{m+1}$
and $1\leq j\leq m+1,$ and thus
\begin{eqnarray*} g_{m+1}(\beta(\Bx))&=&\bigl(g(
\mu_{m}(\pi_{\ne 1} \Bx)),\dots, g(\mu_{m}(\pi_{\ne m+1} \Bx))\bigr)
\\
&\leq_{m+1}&\bigl(\nu_{m}(g_{m}(\pi_{\ne 1} \Bx)),\dots,
\nu_{m}(g_{m}(\pi_{\ne m+1} \Bx)) \bigr )\\
&=&(\nu_{m}(\pi_{\ne 1}g_{m+1}(\Bx)),\dots,
\nu_{m}(\pi_{\ne m+1}g_{m+1}(\Bx)) \bigr)\\
&=&\alpha(g_{m+1}(\Bx))\\
\end{eqnarray*}
and therefore $g_{m+1}\circ \beta\leq_{m+1}\alpha\circ g_{m+1}.$
Since $\nu$ is monotone, $\alpha$ is monotone for $\leq_{m}$
(Theorem \ref{T:mono}), and therefore
\begin{eqnarray*}
g_{m+1}(\beta^{2}(\Bx))&=&(g_{m+1}\beta)(\beta(\Bx))\leq_{m+1}\alpha(g_{m+1}(\beta(\Bx)))\leq
\alpha^{2}(g_{m+1}(\Bx))
\end{eqnarray*}
and inductively
$g_{m+1}(\beta^{n}(\Bx))\leq_{m+1}\alpha^{n}(g_{m+1}(\Bx)))$ for all
$n.$  Since the order is closed,
$$g_{m+1}\bigl(\lim_{n\to\infty}\beta^{n}(\Bx)\bigr)
\leq_{m+1} \lim_{n\to\infty}\alpha^{n}(g_{m+1}(\Bx)).$$ In
particular, $g(\mu_{m+1}(\Bx))\leq \nu_{m+1}(g_{m+1}(\Bx)).$
Induction on $m$ yields the theorem.
\end{proof}

\begin{corollary}\label{C:concav} Let $\mu, \nu$ and $\omega$  be nonexpansive,
coordinatewise contractive $k$-means on complete metric spaces $X,Y,
Z$ respectively.  Let $\leq$ be a closed partial order on $Z$ and
let $g:X\times Y\to Z$ be a continuous function. If $g\circ
(\mu\times \nu)\leq \omega\circ g_{k},$ then $g\circ (\mu\times
\nu)_{n}\leq \omega_{n}\circ g_{n}$ for any $n\geq k$;
\begin{eqnarray*}
&{}&g\Big(\mu_{n}(x_{1},x_{2},\cdots,x_{n}),
\nu_{n}(y_{1},y_{2},\dots,y_{n})\Big)\leq
\omega_{n}\Big(g(x_{1},y_{1}),g(x_{2},y_{2}),\dots,
g(x_{n},y_{n})\Big).
\end{eqnarray*}
\end{corollary}

\begin{proof} By Theorem \ref{T:product}, the product mean
$\mu\times \nu$ is a nonexpansive, coordinatewise contractive
$k$-mean on $X\times Y.$ By Theorem \ref{T:ineq} the inequality
$g\circ (\mu\times \nu)\leq \omega\circ g_{k}$ can be extended to
$n$-means: $g\circ (\mu\times \nu)_{n}\leq \omega_{n}\circ g_{n}.$
From
\begin{eqnarray*}
(\mu\times\nu)_{n}(x_{1},y_{1},x_{2},y_{2},\dots,x_{n},y_{n})&=&\Big(\mu_{n}(x_{1},x_{2},\dots,x_{n}),
\nu_{n}(y_{1},y_{2},\dots,y_{n})\Big)\\
g_{n}(x_{1},y_{1},x_{2},y_{2},\dots,x_{n},y_{n})&=&\Big(g(x_{1},y_{1}),g(x_{2},y_{2}),\dots,
g(x_{n},y_{n})\Big),
\end{eqnarray*}
the proof is completed.
\end{proof}


\section{Means on Hilbert space operators}

In this and the next section we apply and extend our preceding results to the special setting
of positive definite operators on a Hilbert space, in particular to positive definite
Hermitian matrices (in the case the Hilbert space is finite dimensional).

For a Hilbert space $E$, let $\mathcal{B}(E)$ denote the set of
bounded linear operators, $\mathcal{S}(E)\subseteq \mathcal{B}(E)$
the symmetric operators, and $\Omega\subseteq \mathcal{S}(E)$ the
set of positive definite operators on $E$. We define a closed
positive order on $\mathcal{S}(E)$ by $A\leq B$ if $B-A$ is
positive semidefinite. Note that the identity operator $I$ (and
indeed any positive definite operator) is an order unit for
$\mathcal{S}(E)$ (that is, $\mathcal{S}(E)=\bigcup_{n=1}^\infty
[-nI,nI]$, where in general $[A,B]$ denotes the order interval
$[A,B]=\{X\in\mathcal{S}(E): A\leq X\leq B\}$. There is a
corresponding order unit norm given by
$$ \Vert A\Vert=\inf\{t\geq 0: A\in [-tI,tI]\}.$$
This norm generates the same topology on $\mathcal{S}(E)$ as does the operator norm.

We primarily employ the Thompson (or part) metric on $\Omega$ given by
$$d(A,B)=\max\{\log M(A/B),\log M(B/A)\}\mbox{ where } M(A/B):=\inf\{\lambda >0:A\leq
\lambda B\}.$$
A.\ C.\ Thompson \cite{Thomp} has shown that $\Omega$ is a complete metric space with respect to this
metric and the corresponding metric topology on $\Omega$ agrees with the relative
norm topology.  We list some additional elementary properties of the Thompson metric.

\begin{lemma}\label{L:Thompson}
The Thompson metric on the set $\Omega$ of positive definite Hilbert space operators satisfies
\begin{itemize}
\item[(i)] $d(A+B,A+C)\leq d(B,C)$;
\item[(ii)] $A_1\leq A_2$ implies $d(A_1+B,A_1+C)\geq d(A_2+B,A_2+C)$;
\item[(iii)] For $r>0$, $d(rA,rB)=d(A,B)$;
\item[(iv)] $d(A+B,C+D)\leq \max\{d(A,C),d(B,D)\}$;
\item[(v)] $d(A^{-1}, B^{-1})=d(A,B).$
\end{itemize}
\end{lemma}

\begin{proof} (i) There exists $r\geq 1$ such that $\log r=d(B,C)$.  Then
$B\leq rC$, and thus $A+B\leq A+rC\leq rA+rC=r(A+C)$, and
similarly $C\leq rB$ implies $A+C\leq r(A+B)$.  Hence
$d(A+B,A+C)\leq \log r=d(B,C)$.

(ii) The proof of (i) remains valid for $A\geq 0$, and then (ii) follows by rewriting
$A_2$ as $A_1+(A_2-A_1)$.

(iii) For $r>0$, $d(rA,rB)=d(A,B)$ since scalar multiplication by $r$ is an order-isomorphism.

(iv) Suppose that $d(A,C)\leq d(B,D)=\log r.$ Then $B\leq rD, D\leq
rB, A\leq rC, C\leq rA,$ and thus $A+B\leq rC+rD=r(C+D), C+D\leq
rA+rB=r(A+B).$ Hence $d(A+B,C+D)\leq \log r=d(B,D).$

(v) This follows from the order reversing property of operator
inversion.
\end{proof}

The results of Section 6 require forming iterated means from means
that are on the one hand convex, and on the other coordinatewise
$\rho$-contractive and nonexpansive.  The convex mean we employ is
the geometric mean on $\Omega$ defined by
$A\#B=A^{1/2}(A^{-1/2}BA^{-1/2})^{1/2}A^{1/2}$ (see \cite{LL} for
a variety of other characterizations). It is known that the
geometric mean is a convex mean with respect to the Thompson
metric (see, for example, \cite{CPR}, \cite{LL4}).  We record this
fact.

\begin{lemma}\label{L:geometric}
The geometric mean $A\#B$ is a convex mean on $\Omega$ endowed with the respect to the Thompson metric.
\end{lemma}

We close this section by establishing that the arithmetic mean and
the harmonic mean are nonexpansive and coordinatewise
$\rho$-contractive with respect to the Thompson metric when we
restrict to order intervals. This takes some computation.

\begin{lemma}\label{L:arithmetic}
For each $n\in \N$, there exists $\rho_n$, $0<\rho_n< 1$, such that
the arithmetic mean $A(X,Y)=(X+Y)/2$ and the harmonic mean
$H(X,Y)=2(X^{-1}+Y^{-1})^{-1}$ are coordinatewise
$\rho_n$-contractive and nonexpansive on the order interval
$[(1/n)I,nI]\subseteq \Omega$, where $I$ is the identity operator.
\end{lemma}

\begin{proof} By Lemma \ref{L:Thompson} (iv), the arithmetic mean is
nonexpansive. We seek $\rho=\rho_n$ such that for any $A,B,C\in
[(1/n)I,nI]$, $d(A+B,A+C)\leq \rho d(B,C)$ (note that we can drop
the factor of $1/2$ by Lemma \ref{L:Thompson}(iii)).  This is
equivalent to
$$\max\{\log M(A+B/A+C),\log M(A+C/A+B))\}\leq \rho \max\{\log M(B/C),\log M(C/B))\}.$$
Note that if $B\ne C$ (the desired inequality is trivially true if $B=C$),
then either $B\nleq C$ or $C\nleq B$, and then there exists $r>1$ such that
$C\leq rB$, $r=M(C/B)$, and $\log r=d(B,C)$ (or vice-versa with the roles of $B$ and $C$ interchanged).

Suppose now that we can find $\rho$, $0<\rho<1$, such that for all
$A,B,C\in [(1/n)I,nI]$, we have $A+C\leq r^\rho(A+B)$, where
$r=M(C/B)=\max\{M(B/C),M(C/B)\}$.  Then
$$d(A+B,A+C)\leq \log r^\rho=\rho\log r=\rho d(B,C),$$
which is the coordinatewise $\rho$-contractive property.  We
conclude that to establish the coordinatewise $\rho$-contractive
property, it suffices to show the existence of some $\rho$,
$0<\rho<1$, such that for all $A,B, C\in [(1/n)I,nI]$, if $C\leq
rB$ for $1<r\leq n^2$, then $A+C\leq r^\rho(A+B)$. (Note that we
can restrict to $r\leq n^2$ since $C\leq nI=n^2(1/n)I\leq n^2 B$.)
We establish that this is indeed the case by means of the two
following claims.

\noindent\textbf{Claim 1}: \emph{For given $0<\rho<1$ and $r>1$, assume that
$C+(1/n)I\leq r^\rho(B+(1/n)I)$ whenever $C\leq rB$.
Then $A+C\leq r^\rho(A+B)$ for all} $A\geq (1/n)I$.  Indeed
\begin{eqnarray*}
A+C &=& (A-\frac{1}{n}I)+\frac{1}{n}I+C\\
&\leq & (A-\frac{1}{n}I)+r^\rho(\frac{1}{n}I+B)\\
&\leq & r^\rho(A-\frac{1}{n}I)+r^\rho(\frac{1}{n}I+B)=r^\rho(A+B).
\end{eqnarray*}

\noindent\textbf{Claim 2}: \emph{There exists $\rho$, $0<\rho<1$,
such that for all $1<r\leq n^2$, $C+(1/n)I\leq r^\rho(B+(1/n)I)$
whenever} $C\leq rB$, $B,C\in [(1/n)I,nI]$. Indeed for any
$0<\rho< 1$,
$$C+\frac{1}{n}I\leq rB+\frac{1}{n}I=r^\rho B+(r-r^\rho)B+\frac{1}{n}I.$$
For the last two terms we have
$$(r-r^\rho)B+\frac{1}{n}I\leq (r-r^\rho)(nI)+\frac{1}{n}I=
\bigl( (r-r^\rho)n^2+1\bigr)\frac{1}{n}I.$$ To complete the proof of
Claim 2, we need to choose $\rho<1$, but large enough so that
$(r-r^\rho)n^2+1\leq r^\rho$ for $1\leq r\leq n^2$.  The function
$f(x)=x^\rho-n^2(x-x^\rho)-1$ has derivative
$$f'(x)=\rho x^{\rho-1}-n^2(1-\rho x^{\rho-1})=(1+n^2)\rho x^{\rho-1}-n^2.$$
Since the limit of the right-hand expression is $1$ as $\rho\to 1^-$, we conclude that
the derivative is positive for all $x\in [1,n^2]$ for large enough $\rho$ below $1$.
Thus $f$ is increasing on $[1,n^2]$ for $\rho$ near, but below, $1$,
and hence $$(r-r^\rho)n^2+1\leq r^\rho \mbox{ for any } 1\leq r\leq n^2,\ \vert 1-\rho\vert<\ve$$
for some $\ve>0$.

The case of the harmonic mean follows from the fact that the
inversion is an isometry (Lemma \ref{L:Thompson} (v)) and by
applying the preceding result of the arithmetic mean.
\end{proof}

\section{Extending means on Hilbert space operators}

We summarize fundamental results of Kudo and Ando \cite{KA}
(see also \cite{Ba1},  \cite{PT}) concerning operator means and their
relationship to means on the positive reals.  We consider continuous
means on the positive reals, $\mu:\R^+\times \R^+\to \R^+$,
satisfying
\begin{itemize}
\item[(i)]  If $x\leq x'$ and $y\leq y'$, then $\mu(x,y)\leq \mu(x',y')$ (monotonicity);
\item[(ii)] $\mu(tx,ty)=t\mu(x,y)$ for $t,x,y>0$ (homogeneity).
\end{itemize}

A homogeneous two-variable function $\mu$ can be reduced to a
one-variable function $f(x)=\mu(1,x)$.  This reduction defines a
one-to-one correspondence between the continuous means satisfying
(i) and (ii) and the continuous functions $f:\R^+\to\R^+$
satisfying
\begin{itemize}
\item[(O)] $f(1)=1$;
\item[(I)] $f$ is nondecreasing.
\end{itemize}

We consider continuous operator means $\mu$ on $\Omega$, the set
of positive operators on a Hilbert space, satisfying
\begin{itemize}
\item[(a)] If $A\leq A'$, $B\leq B'$, then $\mu(A,B)\leq \mu(A',B')$ (monotonicity);
\item[(b)] $C\mu(A,B)C^*=\mu(CAC^*,CBC^*)$ for all $C$ invertible (the transformer equality),
\end{itemize}
where, as usual, $A\leq B$ means that $B-A$ is positive semidefinite.

The key result of the theory is that the continuous operator means on $\Omega$ satisfying
(a) and (b) are in one-to-one correspondence to the operator monotone
functions $f:\R^+\to\R^+$ satisfying
(O)-(I), where the correspondence $\mu\leftrightarrow f$ is given by
$$\mu(A,B)=A^{1/2}f(A^{-1/2}BA^{-1/2})A^{1/2}.$$
Recall that $f:\R^+\to\R^+$ is \emph{operator monotone} if its
extension to $\Omega$ via the functional calculus is monotone.  It
is this extension that appears in the displayed equality. The
scalar function $f$ is called the \emph{representing function} of
$\mu$.

Basic examples of the preceding theory include (i) the arithmetic
mean on $\R^+$ with representing function $f(x)=(1/2)(1+x)$ and
operator mean $(1/2)(A+B)$, (ii) the geometric mean $\sqrt{ab}$ on
$\R^+$ with representing function $f(x)=\sqrt x$ and operator
geometric mean $A\#B=A^{1/2}(A^{-1/2}BA^{-1/2})^{1/2}A^{1/2}$, and
(iii) the harmonic mean $2/(a^{-1}+b^{-1})$ with representing
function $f(x)=2x/(1+x)$ and operator harmonic mean
$2(A^{-1}+B^{-1})^{-1}$.

We recall from Section 6 and from \cite{KA} the notion of the
iterated composition $\sigma*\tau$ of two operator means $\sigma$
and $\tau$ on $\Omega$. Starting with $\sigma_1=\sigma$ and
$\tau_1=\tau$, we define inductively the sequences of means
$\{\sigma_n\}$ and $\{\tau_n\}$ by $$\sigma_{n+1}(
A,B)=\sigma(\sigma_n(A,B),\tau_n(A,B)),\quad
\tau_{n+1}(A,B)=\tau(\sigma_n(A,b),\tau_n(A,B)).$$ By Theorem 6.2
of \cite{KA} if $\sigma$ and $\tau$ are continuous, monotonic
means satisfying the transformer equality and if at least one is
neither the left- nor the right-trivial mean, then $\{\sigma_n\}$
and $\{\tau_n\}$ converge to the same limiting mean
$\sigma_\infty$, which means that for all $A,B\in \Omega$, $\lim_n
\sigma_n(A,B)=\sigma_\infty(A,B)=\lim_n\tau_n(A,B)$, where the
limit is taken in the weak operator topology.

\begin{theorem}\label{T:compI}
Let $\Omega$ denote the set of positive operators on a Hilbert
space, let $\lambda(A,B)=A\#B$ denote the geometric mean of $A,B$,
and let $\nu$ be continuous, monotonic mean on $\Omega$ satisfying
the transformer equality that is also coordinatewise
$\rho_n$-contractive for $0<\rho_n<1$ and nonexpansive on the order
interval $[(1/n)I,nI]$ for each $n$  with respect to the Thompson
metric.  Then the iterated composition $\mu=\lambda*\nu$ (resp.\
skewed iterated composition) exists, is a coordinatewise
$\rho_n$-contractive, nonexpansive mean when restricted to
$[(1/n)I,nI]$ for each $n$, and hence inductively $\beta$-converges
to a $\beta$-extension for each $n>2$.
\end{theorem}

\begin{proof}
We consider some fixed order interval $\Gamma_n=[(1/n)I,nI]$ and
$A,B\in[(1/n)I,nI]$; note that $[(1/n)I,nI]$ is closed and convex
with respect to any monotonic mean, in particular with respect to
$\lambda$ and $\nu$.  By Proposition \ref{P:approximate} there
exists an iterated composition $\mu_n= \lambda\vert_{\Gamma_n}*
\nu\vert_{\Gamma_n}$ that is nonexpansive and coordinatewise
$\rho_n$-contractive.  Clearly if $m<n$, then $\mu_n$ is an
extension of $\mu_m$.  Thus there exists a unique mean $\mu$ that
extends all of them. Since any $A,B$ belongs to the domain of some
$\mu_n$, $\mu$ is the iterated composition $\lambda*\nu$.

Since $\mu_n$ is nonexpansive and coordinatewise
$\rho_n$-contractive on $\Gamma_n$, it inductively admits a
$\beta$-extension for each $n>2$. It then follows from Corollary
\ref{C:local} that a $\beta$-extension of $\mu$ exists inductively
for each $n>2$.

The case of the skewed iterated composition is analogous.
\end{proof}

The iterated composition of the operator arithmetic mean and
operator geometric mean yields the arithmetic-geometric operator
mean \cite{Fu}. Similarly we have  the harmonic-geometric operator
mean.
\begin{definition}
For two positive definite operators $A,B\in\Omega$ on the Hilbert
space $E$, we define the \emph{arithmetic-geometric mean} or
\emph{Gauss mean} $AGM(A,B)$ to be the iterated composition of the
geometric and arithmetic means, that is, the limit $\lim_n
\lambda_n(A,B)=\lim_n \nu_n(A,B)$, where we define
$\lambda_1(A,B)=A\#B$, the geometric mean, $\nu_1(A,B)=(A+B)/2$,
the arithmetic mean, and inductively
$\lambda_{n+1}(A,B)=\lambda_n(A,B)\#\nu_n(A,B)$ and
$\nu_{n+1}(A,B)=(1/2)(\lambda_n(A,B)+\nu_n(A,B))$.
\end{definition}

It is standard that the preceding iteration defining the
arithmetic-geometric mean of two positive operators converges in
the weak operator topology and agrees with the one arising from
the representing function of arithmetic-geometric mean on the
positive real numbers \cite[Section 6]{KA}.

\begin{corollary} On any order interval $[(1/n)I,nI]$ the
arithmetic-geometric mean $($resp. the harmonic-geometric mean$)$ is
nonexpansive and coordinatewise $\rho_n$-contractive for some
$\rho_n$, $0<\rho_n<1$. Hence on $\Omega$ the arithmetic-geometric
mean $($resp. the harmonic-geometric mean$)$ inductively
$\beta$-extends to an $n$-mean for each $n>2$.
\end{corollary}

\begin{proof} Fix some positive integer $m$.  By Lemma \ref{L:arithmetic}
the arithmetic mean is nonexpansive and coordinatewise
$\rho_m$-contractive on $[(1/m)I,mI]$ for some $0<\rho_m<1$.  We
have already remarked that the geometric mean is a convex mean
with respect to the Thompson metric. The result now follows from
Theorem \ref{T:compI}.
\end{proof}

We briefly recall the operator logarithmic mean as discussed in
\cite{PT}.  The logarithmic mean is defined on $\R^+$ by
$L(a,b)=(b-a)/(\log b-\log a)$. Its representing function is
$f(x)=L(1,x)=(x-1)/\log x$, which is an operator monotone
function.  Hence there exists a corresponding operator logarithmic
mean.  B.\ C.\ Carlson \cite{Ca} has shown that the logarithmic
mean on $\R^+$ is the skewed iterated composition of the geometric
and arithmetic mean.  It then follows from the theory of operator
means as developed in Section 6 of \cite{KA}, particularly Lemma
6.1 and Theorem 6.2, that the operator logarithmic mean is the
corresponding skewed iterated composition of the operator
geometric and arithmetic means on $\Omega$, the set of positive
operators on a Hilbert space, where the limits are taken in the
weak operator topology. However, the same arguments applied in the
previous corollary to the arithmetic-geometric mean viewed as the
iterated composition of the geometric and arithmetic means apply
equally well to logarithmic mean viewed as the skewed iterated
composition of the geometric and arithmetic mean.  We thus obtain
analogously the following

\begin{corollary} On any order interval $[(1/n)I,nI]$ the
logarithmic mean is nonexpansive and coordinatewise
$\rho_n$-contractive for some $\rho_n$, $0<\rho_n<1$.  Hence on
$\Omega$ the logarithmic mean inductively $\beta$-extends to an
$n$-mean for each $n>2$.
\end{corollary}

The preceding corollary provides a positive solution to a problem
raised by Petz and Temesi (\cite{PT}, \cite{Pe}) as to whether the
logarithmic mean $\beta$-converges and hence admits higher
dimensional extensions.

It is easy to obtain  the order relation $L(A,B)\leq AGM(A,B)$
between the logarithmic mean and the arithmetic-geometric mean,
which are monotone (Definition 9.1). Applying Theorem \ref{T:mono}
we have the following

\begin{corollary} The order relation $L_n(A_1,\ldots,A_n)\leq
AGM_n(A_1,\ldots,A_n)$ holds for the extended logarithmic and
arithmetic-geometric $n$-means for each $n>2.$
\end{corollary}

\begin{remark}
In a similar way one can show that $2$-means $\mu$ and $\nu$
that satisfy the inequality $\mu\leq \nu$ and that can both be derived
by some iteration or skew iteration of the arithmetic and  geometric resp.\
harmonic and geometric means satisfy $\mu_n\leq \nu_n$ for all $n>2$.  In this way,
for example, one derives the principal results of \cite{HZY} as a corollary
to our preceding results.
\end{remark}

We remark in closing that a number of ideas in this paper can be
carried over to the study of means on the set of positive elements
of a $C^*$-algebra, particularly by viewing the $C^*$-algebra as a
closed subalgebra of the algebra of bounded operators on a Hilbert
space. For example, one could define the logarithmic mean to be
the skewed iterated composition of the geometric and arithmetic
means and show that it inductively $\beta$-extends to all higher
dimensions.

\end{document}